\documentclass[letterpaper, 10 pt, conference]{IEEEtran}  

\IEEEoverridecommandlockouts                              

\usepackage{epsfig} 
\usepackage{epstopdf}
\usepackage{amsmath} 
\usepackage{amssymb}  
\usepackage{caption}

\usepackage[caption=false,font=footnotesize]{subfig}

\DeclareCaptionLabelSeparator{periodspace}{.\quad}
\title{\LARGE \bf
Global analysis of a geometric PDAV controller by means of coordinate-free linearization}

\author{Michalis Ramp$^{1}$ and Evangelos Papadopoulos$^{2}$ 
\thanks{$^{1}$M. Ramp is with the Department of Mechanical Engineering, National Technical University of Athens, (NTUA) 15780 Athens, Greece
        {\tt\small rampmich@mail.ntua.gr}}%
\thanks{$^{2}$E. Papadopoulos is with the Department of Mechanical Engineering, NTUA, 15780 Athens (tel: +30-210-772-1440; fax: +30-210-772-1455)
        {\tt\small egpapado@central.ntua.gr}}%
}

\begin{document}


\IEEEpubid{Submitted to 2018 European Control Conference}
\IEEEpubidadjcol
\maketitle
\pagestyle{empty}

\begin{abstract}
Tracking a desired Pointing Direction and simultaneously obtaining a reference Angular Velocity (PDAV) around the pointing direction constitutes a very involved and complicated motion encountered in a variaty of robotic, industrial and military applications.
In this paper through the utilization of global analysis and simulation techniques, the smooth closed-loop vector fields induced by the geometric PDAV controller from \cite{VTUAV}, are visualized to gain a deeper understanding of its global stabilization properties.
First through the calculation of a coordinate-free form of the closed-loop linearized dynamics, the local stability of each equilibrium of the system is analyzed.
The results acquired by means of eigenstructure analysis, are used in predicting the frequency of complex precession/nutation oscillations that arise during PDAV trajectory tracking; an important tool in actuator selection.
Finally, by utilizing variational integration schemes, the flow converging to the desired equilibrium and the flow ''close'' to the stable manifold of the saddle equilibrium of the closed-loop system is visualized and analyzed.
Results offer intimate knowledge of the closed-loop vector fields bestowing to the control engineer the ability to anticipate and/or have a rough estimate of the evolution of the solutions.
\end{abstract}

\section{Introduction}
In aerial or underwater robotics, propulsion is obtained frequently by pointing a rotating high-speed propeller in 3D space (vectored actuation) \cite{VTUAV}, \cite{Repoul}.
As the motors of a platform get pointed, they operate in extreme spinning velocities to produce the needed thrust for the platform locomotion, i.e., the motors track a reference Pointing Direction and simultaneously a reference Angular Velocity (PDAV) about the pointing direction.
This PDAV process represents a fundamental control problem for a variety of robotic, industrial and military applications and constitutes a very involved and complicated motion.
Tiltrotor aircrafts and surveillance apparatus like radar or sonar sensors are a few examples, \cite{Osprey}.
In UAVs, this occurs in aerial platforms driven by out-runner motors \cite{VTUAV}.

To obtain an effective PDAV controller to be utilized on general robotic platforms, we studied the PDAV control problem using geometric methods, \cite{Ramp}.
A singularity-free controller was developed, demonstrating improved performance for large initial attitude errors and the ability to negotiate bounded parametric uncertainties.
The modeling and control of an aerial platform (a vectoring tricopter UAV actuated by three out-runner motors that are pointed in 3D space) was studied, with the UAV utilizing the PDAV controller from \cite{Ramp}, along with core modifications to cope with the platforms high precision vectoring requirements \cite{VTUAV}.

The forenamed and a plethora of other geometric control works, study the global closed-loop dynamics of smooth vector fields on nonlinear manifolds, \cite{3DPS}, \cite{3DPI}, \cite{Koditschek}.
Due to the geometric/topological properties of these manifolds, the desired equilibrium in the above systems has an \textit{almost} global domain of attraction that excludes the union of the stable manifolds of its accompanying equilibria \cite{Koditschek}.
By means of global analysis, simulation techniques and 3D visualizations, the influence of these manifolds on the solutions for the attitude control system of a spherical pendulum and a 3D pendulum were investigated, \cite{SMSE}.
The analysis demonstrated the nontrivial influences of those manifolds on the solutions.

\IEEEpubidadjcol
This work is motivated by the need to gain a deeper understanding of the global stabilization properties of the PDAV controller developed in, \cite{VTUAV}, since we intend to employ this controller on an experimental implementation of the aerial platform described in \cite{VTUAV}.
Resultantly the analysis techniques and computational tools described in \cite{SMSE}, are employed in visualizing the smooth closed-loop vector fields of the PDAV controller, in an attempt to obtain an understanding of its global closed-loop properties.
Our investigation shows that the computational approach described in \cite{SMSE} for the visualization of the stable manifold of the saddle equilibrium did not work for the system at hand, but it was effective in producing the flow ''close'' to it, an equally significant result since multiple observations about the system were extracted.
An additional result of this investigation was the development of the capability to estimate the frequency of complex precession/nutation oscillations arising during PDAV trajectory tracking, an important tool for actuator selection.

\section{Kinetics}
\begin{figure}[!h]
\centering
\subfloat{\includegraphics[width=0.7\columnwidth]{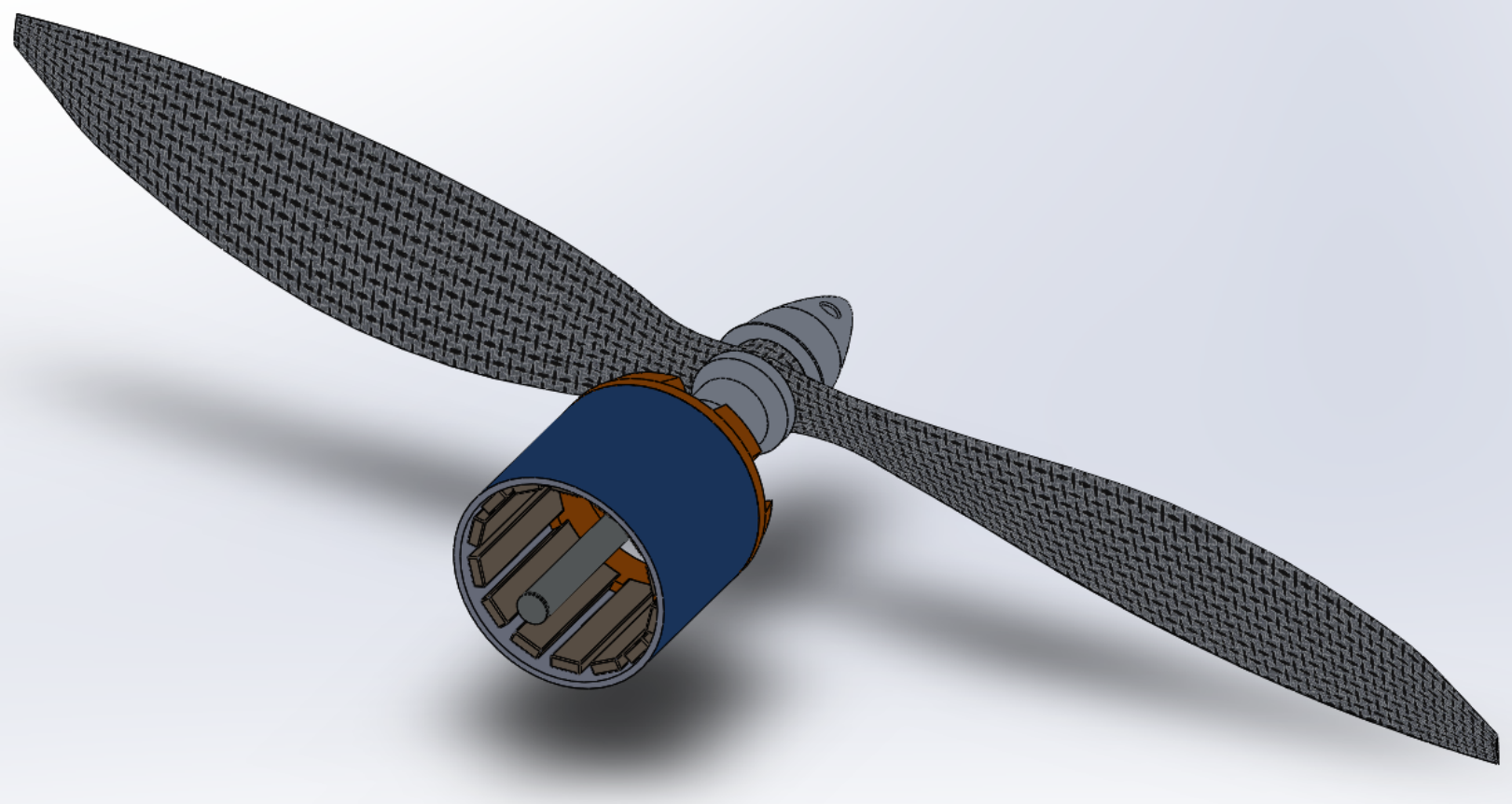}}
\put(-140,80){\parbox{\columnwidth}{$\leftarrow$Propeller}}
\put(-77,52){\parbox{\columnwidth}{$\leftarrow$Propeller holder}}
\put(-162,40){\parbox{\columnwidth}{Flux ring$\rightarrow$}}
\put(-168,20){\parbox{\columnwidth}{Magnets$\longrightarrow$}}

\subfloat{\includegraphics[width=0.9\columnwidth]{./configur-eps-converted-to.pdf}}
\put(-60,80){\parbox{\columnwidth}{${}^{b}\mathbf{u}$}}
\put(-144,50){\parbox{\columnwidth}{$-mg\mathbf{E}_{3}$}}%
\put(-86,47){\parbox{\columnwidth}{$\mathbf{q}$}}
\put(-100,70){\parbox{\columnwidth}{$d\mathbf{R}\mathbf{e}_{3}$}}
\put(-57,59){\parbox{\columnwidth}{$\mathbf{E}_{1}$}}
\put(-148,35){\parbox{\columnwidth}{$\mathbf{E}_{2}$}}
\put(-80,90){\parbox{\columnwidth}{$\mathbf{E}_{3}$}}
\put(-170,42){\parbox{\columnwidth}{$\mathbf{e}_{1}$}}
\put(-138,100){\parbox{\columnwidth}{$\mathbf{e}_{2}$}}
\put(-142,135){\parbox{\columnwidth}{$\mathbf{e}_{3}$}}
\put(-130,121){\parbox{\columnwidth}{${}^{b}\mathbf{F}_{p}$}}
\put(-120,107){\parbox{\columnwidth}{${}^{b}\mathbf{M}_{p}$}}
\caption{
Free body diagram of the out-runner shell/propeller (blue) connected to the inertia frame $\mathbf{I}_{R}$ by a massless axle (green)
}
\label{fig:system}
\end{figure}

The term out-runner motor corresponds to a class of Brushless Direct Current electric motors that spin their outer shell about the stationary windings resulting to a motor that produces far more torque but spin much slower than standard in-runner motors.
We model the moving parts of the out-runner-propeller system i.e., the motor shell with the attached magnets, axle, propeller hub and propeller as a rigid body, see Fig. \ref{fig:system}.
This is done to obtain an understanding of the dynamic phenomena that emerge due to the fast rotations during the pointing procedure.
The system is fully actuated and is shown in Fig. \ref{fig:system}.
It is defined as a rigid body of inertia $\mathbf{J}$ and mass $m$, attached to a frictionless pivot by a massless axle of length $d$ and it is subject to uniform gravity, to the propeller thrust, ${}^{b}\mathbf{F}_{p}$, to the propeller drag torque, ${}^{b}\mathbf{M}_{p}$, and to a control moment ${}^{b}\mathbf{u}\in\mathbb{R}^{3}$ (the components ${}^{b}u_{1}$, ${}^{b}u_{2}$, point the system in 3D space, ${}^{b}u_{3}$, regulates the attitude of the system about the pointing axis i.e., the propeller spinning velocity).
A body fixed frame $\mathbf{I}_{b}\big\{\mathbf{e}_1,\mathbf{e}_2,\mathbf{e}_3\big\}$, attached at the center of mass of the rotating rigid body together with an inertial reference frame $\mathbf{I}_{R}\big\{\mathbf{E}_1,\mathbf{E}_2,\mathbf{E}_3\big\}$, are employed, with $\mathbf{R}(t)\in\text{SO}(3)$ being the rotation matrix from $\mathbf{I}_{b}$ to $\mathbf{I}_{R}$.
The configuration of the system with respect to our task is described by a unit vector $\mathbf{q}(t)\in\mathbb{R}^{3}$,
\begin{IEEEeqnarray}{rCl}
\label{eq:q}
\mathbf{q}(t)&=&\mathbf{R}(t)\mathbf{e}_{3},\mathbf{e}_{3}=[0,0,1]^{T}\label{eq:q}\IEEEeqnarraynumspace
\end{IEEEeqnarray}
collinear with the axis of body rotation, and by the component of the angular velocity about $\mathbf{q}(t)$ given by ${}^{b}\boldsymbol{\omega}(t)\cdot\mathbf{e}_{3}\in\mathbb{R}$.

The configuration space can be either described by an element of $\text{S}^2=\{ \mathbf{q}\in\mathbb{R}^{3}\lvert\mathbf{q}^{T}\mathbf{q}=1 \}$ (the two-sphere) or by an element of the special orthogonal group $\text{SO}(3)=\{ \mathbf{R}\in\mathbb{R}^{3\times 3}\lvert\mathbf{R}^{T}\mathbf{R}=\mathbf{I},\text{det}[\mathbf{R}]=1\}$, even though the attitude about the pointing direction is irrelevant.
The plane tangent to the unit sphere at $\mathbf{q}$ is the tangent space $\text{T}_q\text{S}^2{=}\{ \boldsymbol{\xi}{\in}\mathbb{R}^{3}\lvert\mathbf{q}^{T}\boldsymbol{\xi}=0 \}$.

The equations of motion of the attitude dynamics are,
\begin{IEEEeqnarray}{rCl}
\label{eq:eom}
\mathbf{J}{}^{b}\dot{\boldsymbol{\omega}}+S({}^{b}\boldsymbol{\omega})\mathbf{J}{}^{b}\boldsymbol{\omega}&=&{}^{b}\mathbf{u}+{}^{b}\mathbf{M}_{p}-S(d\mathbf{e}_{3})mg\mathbf{R}^{T}\mathbf{E}_{3}\IEEEyesnumber\IEEEyessubnumber\label{eq:eom_1}\\
\dot{\mathbf{R}}&=&\mathbf{R}S({}^{b}\boldsymbol{\omega})\IEEEyessubnumber\label{eq:eom_2}
\end{IEEEeqnarray}
and the constant $g$ is the gravitational acceleration.
The cross product map, $S(.):\mathbb{R}^{3}\rightarrow\mathfrak{so}(3)$, and its inverse map, $S^{-1}(.):\mathfrak{so}(3)\rightarrow\mathbb{R}^{3}$, are defined in the Appendix.
Using (\ref{eq:q}) and (\ref{eq:eom_2}), the rate of change of $\mathbf{q}$ is given by,
\begin{IEEEeqnarray}{C}
\dot{\mathbf{q}}=S(\mathbf{R}{}^{b}\boldsymbol{\omega})\mathbf{q}=\mathbf{R}S({}^{b}\boldsymbol{\omega})\mathbf{e}_{3}\label{eq:dq}
\end{IEEEeqnarray}

\section{Control System}

Before experimentally applying the PDAV controller \cite{VTUAV}, deep understanding of its closed-loop properties must be obtained; to this end, the controller is summarized next.
For a thorough derivation of the controller see \cite{VTUAV}, \cite{Ramp}.
The error function, \cite{RBAC},
\begin{IEEEeqnarray}{c}
\Psi(\mathbf{q},\mathbf{q}_{d})=1-\mathbf{q}^{T}\mathbf{q}_{d}\label{eq:psi}
\end{IEEEeqnarray}
yields the attitude and angular velocity error vectors, \cite{VTUAV},
\begin{IEEEeqnarray}{rCl}
{}^{b}\mathbf{e}_{q}(\mathbf{R},\mathbf{R}_{d})&=&\mathbf{R}^{T}S(\mathbf{q}_{d})\mathbf{q}\label{eq:eq}\\
^{b}\mathbf{e}_{\omega}({}^{b}\boldsymbol{\omega},{}^{b}\boldsymbol{\omega}_{d},\mathbf{R},\mathbf{R}_{d})&=&{}^{b}\boldsymbol{\omega}-\mathbf{R}^{T}\mathbf{R}_{d}{}^{b}\boldsymbol{\omega}_{d}\label{eq:eo}
\end{IEEEeqnarray}
The control law for a desired pointing direction $\mathbf{q}_{d}=\mathbf{R}_{d}\mathbf{e}_{3}\in\text{S}^{2}$ and a desired angular velocity ${}^{b}\boldsymbol{\omega}_{d}=\omega_{d}\mathbf{e}_{3}$, is given by,
\begin{IEEEeqnarray}{rCl}
\label{eq:u}
{}^{b}\mathbf{u}&=&\eta^{-1}\mathbf{J}\left(-\eta\boldsymbol{\alpha}{-}(\Lambda{+}\Psi){}^{b}\dot{\mathbf{e}}_{q}{-}\dot{\Psi}{}^{b}\mathbf{e}_{q}{-}\gamma\mathbf{s}\right){-}{}^{b}\mathbf{f}\IEEEyesnumber\IEEEyessubnumber\label{eq:u_1}\\
\boldsymbol{\alpha}&=&S({}^{b}\boldsymbol{\omega})\mathbf{R}^{T}\mathbf{R}_{d}{}^{b}\boldsymbol{\omega}_{d}-\mathbf{R}^{T}\mathbf{R}_{d}{}^{b}\dot{\boldsymbol{\omega}}_{d}\IEEEyessubnumber\\
{}^{b}\mathbf{f}&=&{}^{b}\mathbf{M}_{p}-S(d\mathbf{e}_{3})mg\mathbf{R}^{T}\mathbf{E}_{3}-S({}^{b}\boldsymbol{\omega})\mathbf{J}{}^{b}\boldsymbol{\omega}\IEEEyessubnumber\\
\mathbf{s}&=&(\Lambda+\Psi){}^{b}\mathbf{e}_{q}+\eta{}^{b}\mathbf{e}_{\boldsymbol{\omega}}\IEEEyessubnumber\label{eq:S}
\end{IEEEeqnarray}
where $\eta,\gamma,\Lambda\in\mathbb{R}^{+}$, $\omega_{d}\in\mathbb{R}$, and the terms $\dot{\Psi}$ and ${}^{b}\dot{\mathbf{e}}_{q}$ are given in the Appendix by (\ref{eq:dpsi}), (\ref{eq:deq}).

The closed loop dynamics under the action of (\ref{eq:u}) are,
\begin{IEEEeqnarray}{rCl}
\label{eq:cld}
{}^{b}\dot{\boldsymbol{\omega}}&=&\eta^{-1}\left({-}(\Lambda{+}\Psi){}^{b}\dot{\mathbf{e}}_{q}{-}\dot{\Psi}{}^{b}\mathbf{e}_{q}{-}\gamma\mathbf{s}\right)-\boldsymbol{\alpha}\IEEEyesnumber\IEEEyessubnumber\label{eq:cld_1}\\
\dot{\mathbf{q}}&=&\mathbf{R}S({}^{b}\boldsymbol{\omega})\mathbf{e}_{3}\IEEEyessubnumber\label{eq:cld_2}
\end{IEEEeqnarray}
Since $\mathbf{q}_{d}{=}\mathbf{R}_{d}\mathbf{e}_{3},{}^{b}\boldsymbol{\omega}_{d}{=}\omega_{d}\mathbf{e}_{3},{}^{b}\dot{\boldsymbol{\omega}}_{d}{=}\mathbf{0}$, the set of admissible closed loop equilibria solutions of (\ref{eq:cld}), is given by,
\begin{IEEEeqnarray}{rCl}
\label{eq:cle}
\!\!\!(\mathbf{q}_{e}{,}{}^{b}\boldsymbol{\omega}_{e})&{\in}&\{{\big(}{-}\mathbf{q}_{d}{=}{\exp}(\pi{S(\mathbf{R}_{d}\mathbf{e}_{1})}){\exp}(\zeta S(\mathbf{R}_{d}\mathbf{e}_{3}))\mathbf{R}_{d}\mathbf{e}_{3}\IEEEnonumber\\
&&{,}\!{-}{}^{b}\boldsymbol{\omega}_{d}{\big)}\!{,}\!{\big(}\mathbf{q}_{d}{=}{\exp}(\zeta S(\mathbf{R}_{d}\mathbf{e}_{3}))\mathbf{R}_{d}\mathbf{e}_{3}{,}\!{}^{b}\boldsymbol{\omega}_{d}{\big)}\!{|}\zeta{\in}\mathbb{R}\}
\end{IEEEeqnarray}
with the first element of (\ref{eq:cle}) to correspond to the antipodal equilibrium and the second to the desired equilibrium.

Finally using the Lyapunov function,
\begin{IEEEeqnarray}{C}
V(\Psi,{}^{b}\mathbf{e}_{q},{}^{b}\mathbf{e}_{\omega})=\frac{1}{2}\mathbf{s}^{T}\mathbf{s}=\frac{1}{2}\lVert\mathbf{s}\rVert^{2}\label{eq:V}
\end{IEEEeqnarray}
it was shown that the desired equilibrium $(\mathbf{q}_{d},{}^{b}\boldsymbol{\omega}_{d})$ is almost globally exponentially stable, \cite{VTUAV}, \cite{Ramp}.

\section{A Pdav Tracking Case\label{sec:pdav}}

To underline the importance of modeling the system as in Fig. \ref{fig:system} and to showcase the rich dynamic phenomena that arise during pointing an out-runner motor, a simulation is first presented.
The system is initially at equilibrium i.e., $(\mathbf{q}(0){=}\mathbf{e}_{3},{}^{b}\boldsymbol{\omega}{=}1000\mathbf{e}_{3}rad/s)$.
The trajectory performed under the action of (\ref{eq:u}) is that of pointing $90^{o}$ about $\mathbf{E}_{1}$ axis, $90^{o}$ about $\mathbf{E}_{3}$ axis, while the spinning velocity is maintained at $1000rad/s$. 
The trajectories are produced using ''minimum snap'' polynomials \cite{Mellinger}.
The gains $\Lambda,\eta,\gamma$ of (\ref{eq:u}) are,
\begin{IEEEeqnarray}{rCl}
\Lambda=25{\cdot}10^{6},\:\eta=12{\cdot}10^{3},\:\gamma=500\label{eq:gains}
\end{IEEEeqnarray}
The inertial matrix used was obtained by a CAD design of the out-runner shell/propeller assembly of Fig. \ref{fig:system} and is given by,
\begin{IEEEeqnarray}{L}
\mathbf{J}=\begin{bmatrix}
3.612&0.762&0\\
0.762&8.709&0\\
0&0&6.076
\end{bmatrix}{\cdot}10^{-5}[kgm^2]
\label{eq:Ji}
\end{IEEEeqnarray}

\begin{figure}[!h]
\centering
\subfloat[\label{fig:pdav_1}]{\includegraphics[width=0.45\columnwidth]{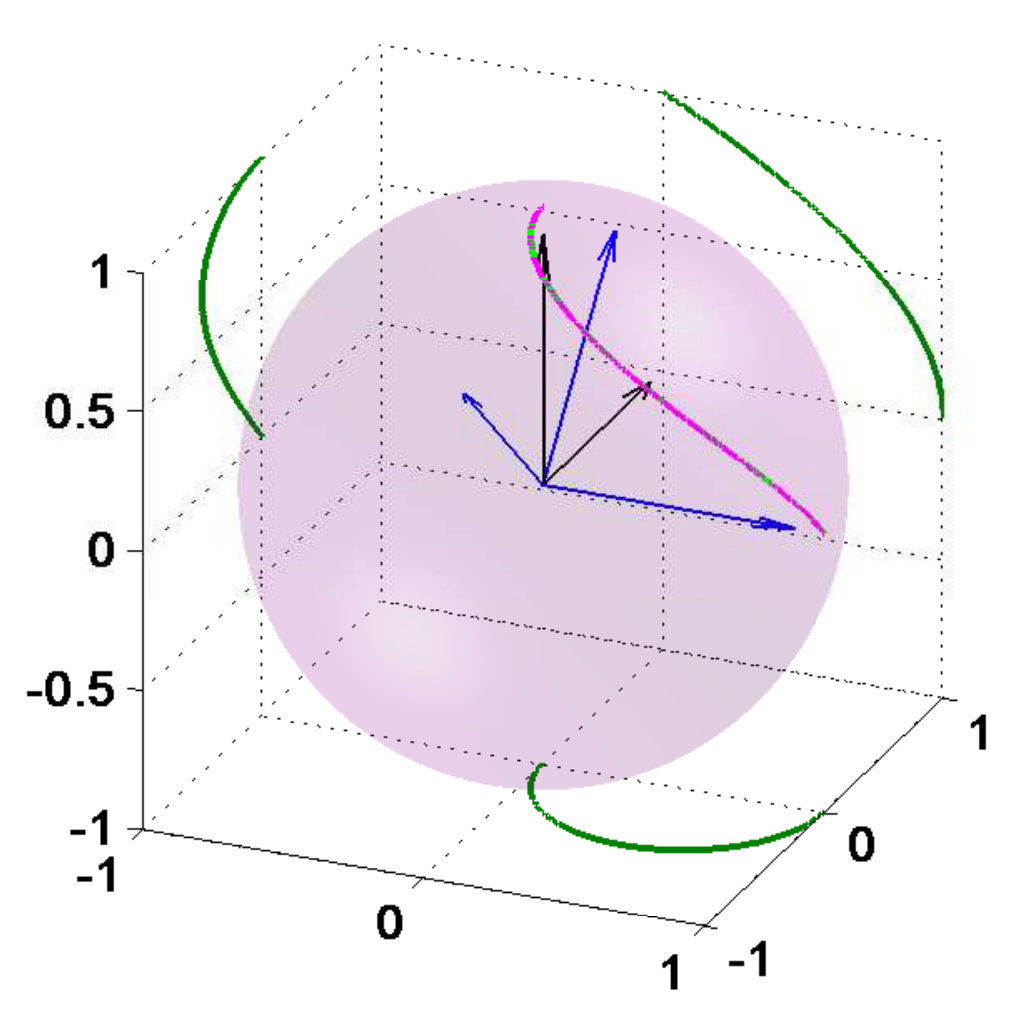}}
\put(-20,55){\parbox{0.49\columnwidth}{$\mathbf{E}_{1}$}}
\put(-40,70){\parbox{0.49\columnwidth}{$\mathbf{E}_{2}$}}
\put(-58,93){\parbox{0.49\columnwidth}{$\mathbf{E}_{3}$}}
\put(-30,48){\parbox{0.49\columnwidth}{$\mathbf{e}_{3}$}}
\put(-45,88){\parbox{0.49\columnwidth}{$\mathbf{e}_{1}$}}
\put(-72,72){\parbox{0.49\columnwidth}{$\mathbf{e}_{2}$}}
~
\subfloat[\label{fig:pdav_2}]{\includegraphics[width=0.49\columnwidth]{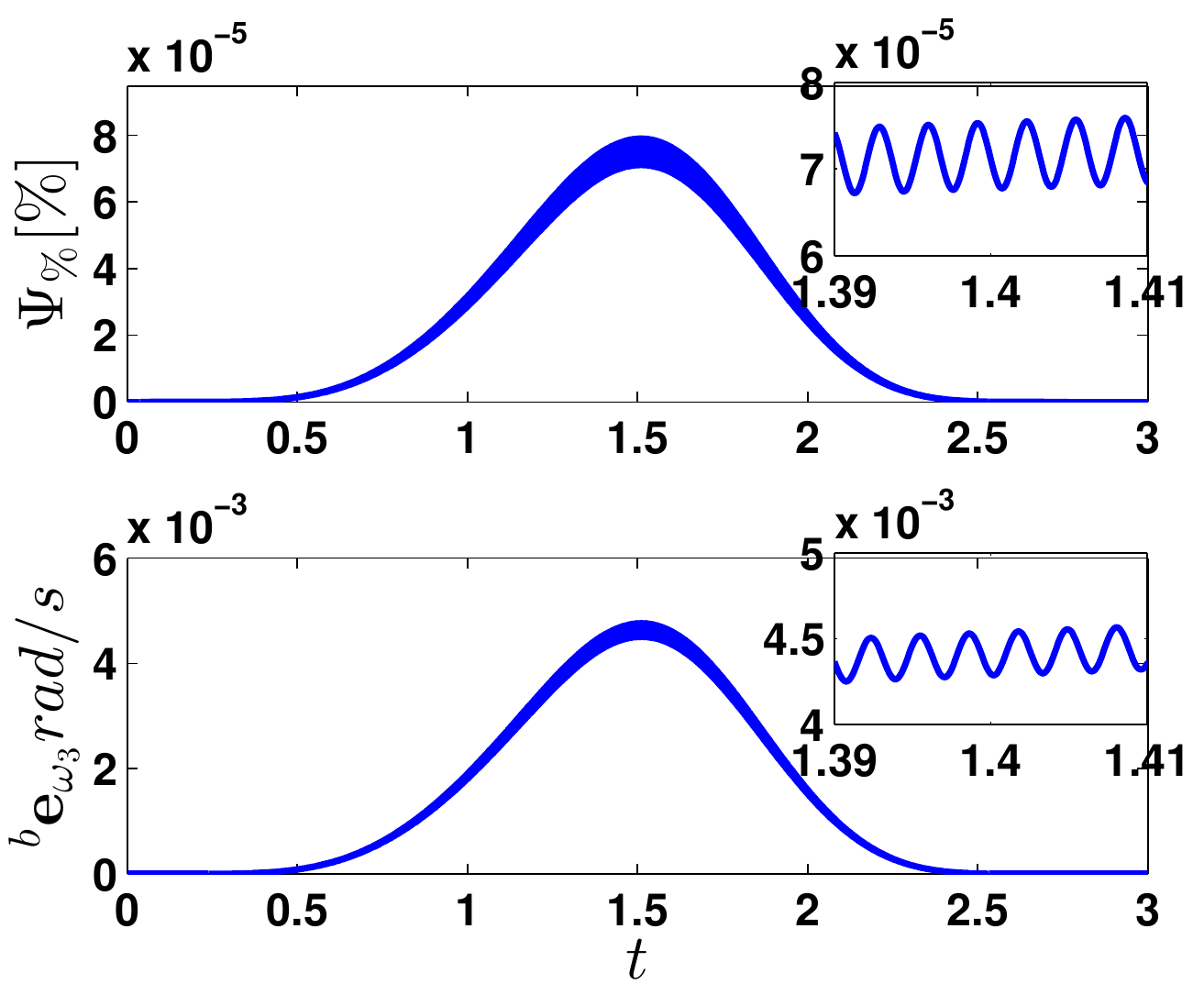}}

\subfloat[\label{fig:pdav_3}]{\includegraphics[width=0.49\columnwidth]{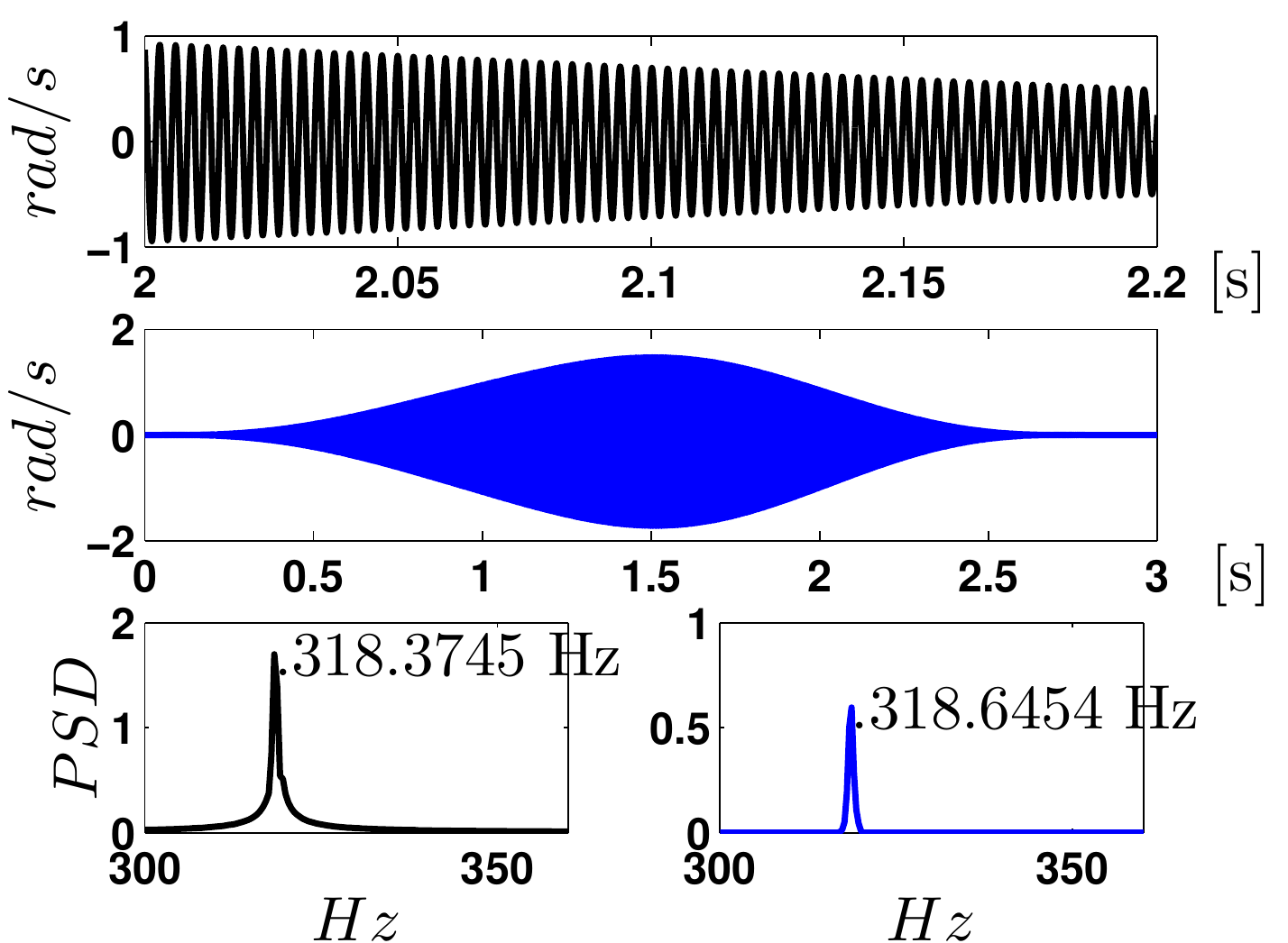}}
~
\subfloat[\label{fig:pdav_4}]{\includegraphics[width=0.49\columnwidth]{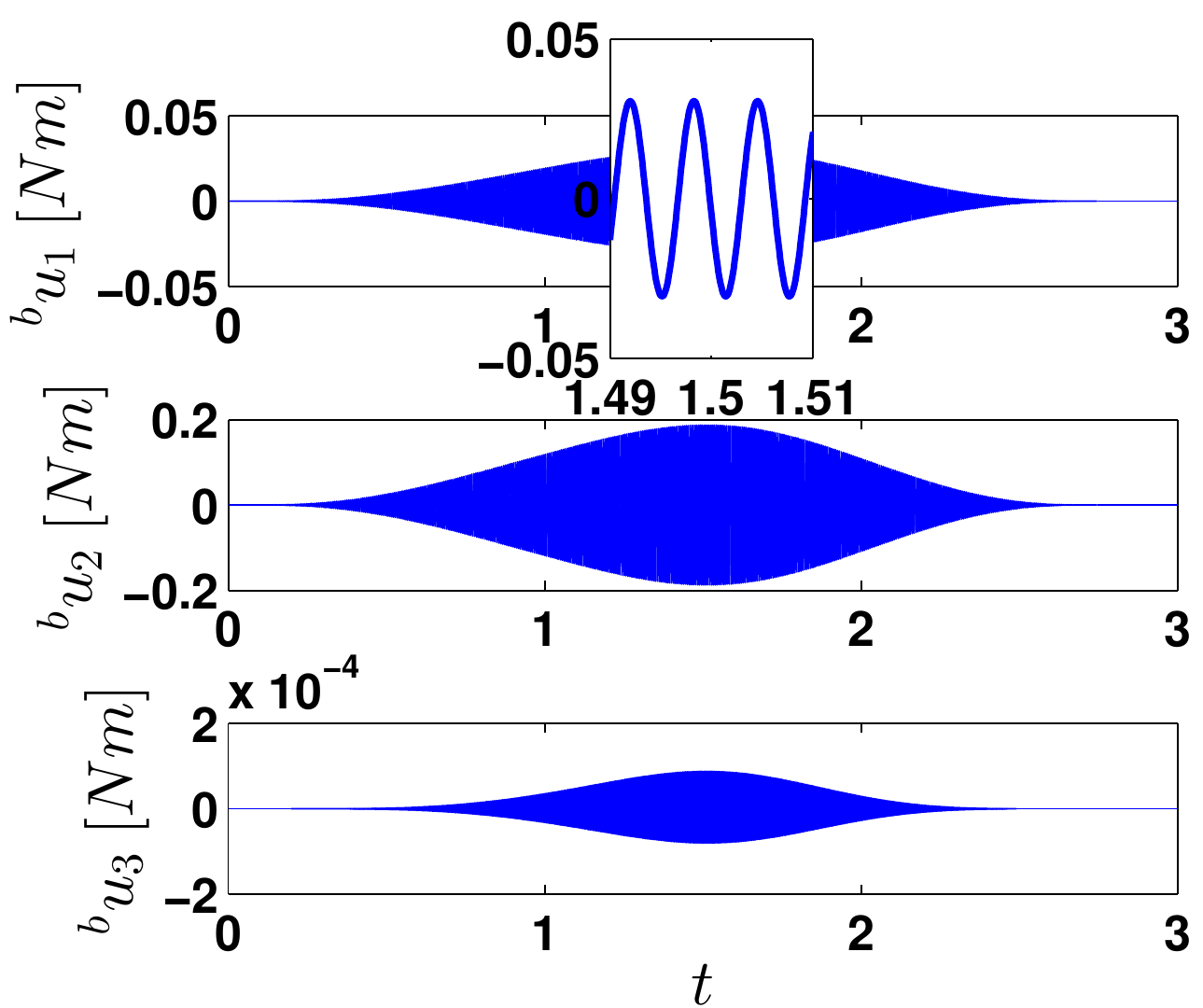}}
\caption{PDAV trajectory with steady propeller velocity at 1000$rad/s$.
(\ref{fig:pdav_1}) Attitude maneuver with projections (green).
(\ref{fig:pdav_2}) Top: Percentage pointing error by (\ref{eq:psi}), Bottom: Propeller spin error. 
(\ref{fig:pdav_3}) Top: Precession rate $\dot{\phi}$ (black). Middle: Nutation rate $\dot{\theta}$ (blue). Bottom Left: Precession frequency (black) by FFT. Bottom Right: Nutation frequency (blue) by FFT. 
(\ref{fig:pdav_4}) Control effort.
}
\label{fig:pdav}
\end{figure}

The closed-loop 3D pointing response (including the initial/final attitude) is shown in Fig. \ref{fig:pdav_1}, showcasing a smooth maneuver.
The controller (\ref{eq:u}) achieves high precision tracking of the desired pointing direction and propeller spinning velocity as indicated in Fig. \ref{fig:pdav_2}.
The percentage attitude error (using (\ref{eq:psi})) remains below $\Psi_{\%}{=}8{\cdot}10^{-5}\%$, see Fig. \ref{fig:pdav_2} (top row) (note that $\Psi_{\%}{=}100\%$ at the maximum pointing error corresponding to $180^{o}$ wrt., an axis angle rotation i.e., the antipodal equilibrium).
The desired propeller speed is tracked faithfully with the corresponding tracking error component, ${}^{b}\mathbf{e}_{\omega_{3}}$, to remain below $6{\cdot}10^{-3}rad/s$, see Fig. \ref{fig:pdav_2} (bottom row).

Of significant importance is the emergence of non trivial nutation/precession dynamic oscillations, shown in Fig. \ref{fig:pdav_3} (top and middle rows).
\textit{Nutation} i.e., a swaying/nodding motion of the pointing axis, is described by a change of the second Euler angle, $\theta$, of the ''313'' sequence.
\textit{Precession} i.e., the azimuth variation of the pointing axis about $\mathbf{E}_{3}$, is expressed by a change of the first Euler angle, $\phi$.
An FFT analysis conducted on the nutation/precession signals reveals high frequency oscillations of 318.6454 Hz (see Fig. \ref{fig:pdav_3} (bottom row)) that the controller compensates for continuously (see Fig. \ref{fig:pdav_4}) to achieve the smooth pointing response shown in Fig. \ref{fig:pdav_1}.
Moreover the torque requirements are revealed as the control effort in Fig. \ref{fig:pdav_4} shows that pointing actuators must generate up to 0.2Nm despite the small inertia of the assembly.
In the majority of the literature involving pointing or tilting out-runner motors, the motor is modeled simply as a source of thrust and torque.
As a result this dynamically rich response is lost, an important omission, since these phenomena play an important role during experimental implementations.

\section{Linearization}

In the previous section, the importance of proper modeling was shown.
The goal of this work is to study the closed-loop equilibrium properties explicitly; thus a coordinate-free form of the linearized dynamics of (\ref{eq:cld}) is developed and the local stability of each equilibrium is analyzed.
To this end, the closed-loop equations are linearized about each equilibrium using suitable expressions for the variation of the states, ensuring that the perturbation of the equilibrium lies on the configuration space.
This is achieved as in \cite{RBAC}, by using the exponential map, (\ref{expm}), to define the perturbation of the equilibrium in terms of a perturbation parameter $\epsilon$ as,
\begin{IEEEeqnarray}{rCl}
\mathbf{q}(t,\epsilon)&=&\exp(\epsilon S(\boldsymbol{\xi}))\mathbf{q}(t)\label{eq:del_q}\\
\mathbf{R}(t,\epsilon)&=&\exp(\epsilon S(\boldsymbol{\xi}))\mathbf{R}(t)\label{eq:delR}
\end{IEEEeqnarray}
with $\boldsymbol{\xi}\in\text{T}_q\text{S}^2$.
The perturbation of the angular velocity is,
\begin{IEEEeqnarray}{rCl}
{}^{b}\boldsymbol{\omega}(t,\epsilon)={}^{b}\boldsymbol{\omega}(t)+\epsilon\delta\mathbf{w}(t)\label{eq:delo}
\end{IEEEeqnarray}
and the curve $\delta\mathbf{w}(t)\in\mathbb{R}^{3}$.
Note that if the perturbation parameter $\epsilon=0$ then $(\mathbf{q}(0,0),{}^{b}\boldsymbol{\omega}(0,0))=(\mathbf{q}_{e},{}^{b}\boldsymbol{\omega}_{e})$, meaning that, $(\mathbf{q}(t,0),{}^{b}\boldsymbol{\omega}(t,0))=(\mathbf{q}_{e},{}^{b}\boldsymbol{\omega}_{e}),\forall t\in[0,\infty)$.
To save space, onwards the dependency on time $t$ is dropped.
By utilizing (\ref{eq:del_q}) the infinitesimal variation of $\mathbf{q}(t,\epsilon)$ is,
\begin{IEEEeqnarray}{c}
\delta \mathbf{q}=\frac{d}{d\epsilon}\Big\lvert_{\epsilon=0}\exp(\epsilon S(\boldsymbol{\xi}))\mathbf{q}=S(\boldsymbol{\xi})\mathbf{q}\label{eq:delq}
\end{IEEEeqnarray}

To obtain the coordinate-free form of the closed-loop dynamics we substitute (\ref{eq:del_q})-(\ref{eq:delo}) into (\ref{eq:cld}), differentiate both sides of the resulting equation wrt., $\epsilon$, and substitute $\epsilon=0$.

First by differentiating (\ref{eq:delq}) we get,
\begin{IEEEeqnarray}{rCl}
\delta\dot{\mathbf{q}}=S(\dot{\boldsymbol{\xi}})\mathbf{q}+S(\boldsymbol{\xi})\dot{\mathbf{q}}\label{eq:Ddelq}
\end{IEEEeqnarray}
Substituting (\ref{eq:del_q})-(\ref{eq:delo}) in (\ref{eq:cld_2}) and differentiating wrt., $\epsilon$,
\begin{IEEEeqnarray}{rCl}
\delta\dot{\mathbf{q}}=S(\boldsymbol{\xi})(S(\mathbf{R}{}^{b}\boldsymbol{\omega})\mathbf{q})+S(\mathbf{R}\delta\mathbf{w})\mathbf{q}\label{eq:Ddelqalt}
\end{IEEEeqnarray}
Equating (\ref{eq:Ddelq}) with (\ref{eq:Ddelqalt}), and substituting (\ref{eq:cld_2}) we get,
\begin{IEEEeqnarray}{rCl}
S(\dot{\boldsymbol{\xi}})\mathbf{q}=S(\mathbf{R}\delta\mathbf{w})\mathbf{q}\label{eq:Ddelqint_1}
\end{IEEEeqnarray}
To be able to solve (\ref{eq:Ddelqint_1}) wrt., $\dot{\boldsymbol{\xi}}$ since both sides are perpendicular to $\mathbf{q}$ it holds that, $S(\mathbf{q})(S(\dot{\boldsymbol{\xi}})\mathbf{q})=S(\mathbf{q})(S(\mathbf{R}\delta\mathbf{w})\mathbf{q})$.
Using the identity $S(\mathbf{q})(S(\dot{\boldsymbol{\xi}})\mathbf{q})=(\mathbf{q}^{T}\mathbf{q})\dot{\boldsymbol{\xi}}-(\mathbf{q}^{T}\dot{\boldsymbol{\xi}})\mathbf{q}$ then,
\begin{IEEEeqnarray}{rCl}
\dot{\boldsymbol{\xi}}-(\mathbf{q}^{T}\dot{\boldsymbol{\xi}})\mathbf{q}=S(\mathbf{q})(S(\mathbf{R}\delta\mathbf{w})\mathbf{q})\label{eq:Ddelqint_2}
\end{IEEEeqnarray}
Since $\boldsymbol{\xi}\in\text{T}_q\text{S}^2$ then $\dot{\boldsymbol{\xi}}^{T}\mathbf{q}=-\boldsymbol{\xi}^{T}\dot{\mathbf{q}}$ and by applying it to (\ref{eq:Ddelqint_2}) and rearranging terms the linearized equation of motion of (\ref{eq:cld_2}) is obtained as,
\begin{IEEEeqnarray}{rCl}
\dot{\boldsymbol{\xi}}=\mathbf{q}\mathbf{q}^{T}S(\mathbf{R}{}^{b}\boldsymbol{\omega})\boldsymbol{\xi}+(\mathbf{I}-\mathbf{q}\mathbf{q}^{T})\mathbf{R}\delta^{b}\boldsymbol{\omega}\label{eq:Ddelqint_3}
\end{IEEEeqnarray}

In the same manner utilizing (\ref{eq:del_q})-(\ref{eq:delo}) in (\ref{eq:cld_1}) and differentiating wrt., $\epsilon$, after considerable manipulations, we get the linearized equation of motion of (\ref{eq:cld_1}) as,
\begin{IEEEeqnarray}{rCl}
\delta^{b}\dot{\boldsymbol{\omega}}&=&\mathbf{R}^{T}S(\mathbf{R}_{d}{}^{b}\dot{\boldsymbol{\omega}}_{d})\boldsymbol{\xi}\IEEEnonumber\\
&&-\Big(S({}^{b}\boldsymbol{\omega})\mathbf{R}^{T}S(\mathbf{R}_{d}{}^{b}\boldsymbol{\omega}_{d})\boldsymbol{\xi}-S(\mathbf{R}^{T}\mathbf{R}_{d}{}^{b}\boldsymbol{\omega}_{d})\delta^{b}\boldsymbol{\omega}\Big)\IEEEnonumber\\
&&-\eta^{-1}\Big({}^{b}\dot{\mathbf{e}}_{q}\mathbf{q}_{d}^{T}S(\mathbf{q})\boldsymbol{\xi}+(\Lambda+\Psi){}^{b}\dot{\mathbf{e}}_{q}^{\epsilon}\Big)\IEEEnonumber\\
&&-\eta^{-1}\Big({}^{b}\mathbf{e}_{q}(\mathbf{R}{}^{b}\mathbf{e}_{\omega})^{T}S(\mathbf{R}{}^{b}\mathbf{e}_{q})\boldsymbol{\xi}+{}^{b}\mathbf{e}_{q}(\mathbf{R}{}^{b}\mathbf{e}_{\omega})^{T}\mathbf{R}{}^{b}\mathbf{e}_{q}^{\epsilon}\IEEEnonumber\\
&&+{}^{b}\mathbf{e}_{q}(\mathbf{R}{}^{b}\mathbf{e}_{q})^{T}S(\mathbf{R}{}^{b}\mathbf{e}_{\omega})\boldsymbol{\xi}{+}{}^{b}\mathbf{e}_{q}(\mathbf{R}{}^{b}\mathbf{e}_{q})^{T}\mathbf{R}{}^{b}\mathbf{e}_{\omega}^{\epsilon}{+}\dot{\Psi}{}^{b}\mathbf{e}_{q}^{\epsilon}\Big)\IEEEnonumber\\
&&-\gamma\eta^{-1}\Big({}^{b}\mathbf{e}_{q}\mathbf{q}^{T}_{d}S(\mathbf{q})\boldsymbol{\xi}+(\Lambda+\Psi){}^{b}\mathbf{e}_{q}^{\epsilon}+\eta{}^{b}\mathbf{e}_{\omega}^{\epsilon}\Big)\label{eq:Ddelo}
\end{IEEEeqnarray}
and ${}^{b}\mathbf{e}_{q}^{\epsilon},{}^{b}\dot{\mathbf{e}}_{q}^{\epsilon},{}^{b}\mathbf{e}_{\omega}^{\epsilon}$ are in the Appendix as (\ref{eq:deqe}), (\ref{eq:Ddeqe}) and (\ref{eq:deoe}).
Finally, because (\ref{eq:Ddelqint_3}), (\ref{eq:Ddelo}) will be used to extract the eigenvalues at the equilibria they are further rearranged to get,
\begin{IEEEeqnarray}{rCl}
\dot{\mathbf{x}}=\begin{bmatrix}
\dot{\boldsymbol{\xi}}\\
\delta^{b}\dot{\boldsymbol{\omega}}
\end{bmatrix}&=&
\begin{bmatrix}
\mathbf{\Xi}_{\xi}&\mathbf{\Xi}_{\omega}\\
\mathbf{\Omega}_{\xi}&\mathbf{\Omega}_{\omega}
\end{bmatrix}
\begin{bmatrix}
\boldsymbol{\xi}\\
\delta^{b}\boldsymbol{\omega}
\end{bmatrix}=\mathbf{A}\mathbf{x}\label{eq:mat}
\end{IEEEeqnarray}
and the terms of the matrix $\mathbf{A}\in\mathbb{R}^{6\times6}$ are given by,

\begin{IEEEeqnarray*}{rCl}
\mathbf{\Xi}_{\xi}&=&\mathbf{q}\mathbf{q}^{T}S(\mathbf{R}{}^{b}\boldsymbol{\omega}),\mathbf{\Xi}_{\omega}=(\mathbf{I}-\mathbf{q}\mathbf{q}^{T})\mathbf{R}\\
\mathbf{\Omega}_{\xi}&=&\mathbf{R}^{T}S(\mathbf{R}_{d}{}^{b}\dot{\boldsymbol{\omega}}_{d}){-}S({}^{b}\boldsymbol{\omega})\mathbf{R}^{T}S(\mathbf{R}_{d}{}^{b}\boldsymbol{\omega}_{d})\\
&&{-}{\eta}^{-1}\big\{{}^{b}\mathbf{e}_{q}(\mathbf{R}{}^{b}\mathbf{e}_{\omega})^{T}S(\mathbf{R}{}^{b}\mathbf{e}_{q})\\
&&+{}^{b}\mathbf{e}_{q}(\mathbf{R}{}^{b}\mathbf{e}_{q})^{T}S(\mathbf{R}{}^{b}\mathbf{e}_{\omega}){+}{}^{b}\dot{\mathbf{e}}_{q}\mathbf{q}_{d}^{T}S(\mathbf{q})\\
&&{+}{(\Lambda{+}\Psi)}\Big(\mathbf{R}^{T}S(S(\dot{\mathbf{q}}_{d})\mathbf{q}{+}S(\mathbf{q}_{d})\dot{\mathbf{q}})
\\
&&{-}\mathbf{R}^{T}S(\mathbf{q}_{d})S(\mathbf{R}S({}^{b}\boldsymbol{\omega})\mathbf{e}_{3})
{-}\mathbf{R}^{T}S(\dot{\mathbf{q}}_{d})S(\mathbf{q})\\
&&-S({}^{b}\boldsymbol{\omega})\mathbf{R}^{T}\big(S(S(\mathbf{q}_{d})\mathbf{q})-S(\mathbf{q}_{d})S(\mathbf{q})\big)\Big)\\
&&{+}{\left({}^{b}\mathbf{e}_{q}(\mathbf{R}{}^{b}\mathbf{e}_{\omega})^{T}\mathbf{R}{+}{\left(\dot{\Psi}{+}\gamma(\Lambda{+}\Psi)\right)}\mathbf{I}\right)}\cdot\\
&&\mathbf{R}^{T}\Big(S(S(\mathbf{q}_{d})\mathbf{q})-S(\mathbf{q}_{d})S(\mathbf{q})\Big){+}{\gamma}{}^{b}\mathbf{e}_{q}\mathbf{q}^{T}_{d}S(\mathbf{q})\\
&&{-}({}^{b}\mathbf{e}_{q}(\mathbf{R}{}^{b}\mathbf{e}_{q})^{T}\mathbf{R}{+}\eta\gamma\mathbf{I})\mathbf{R}^{T}S(\mathbf{R}_{d}{}^{b}\boldsymbol{\omega}_{d})\big\}\\
\mathbf{\Omega}_{\omega}&=&S(\mathbf{R}^{T}\mathbf{R}_{d}{}^{b}\boldsymbol{\omega}_{d}){-}{\eta}^{-1}\big\{({}^{b}\mathbf{e}_{q}(\mathbf{R}{}^{b}\mathbf{e}_{q})^{T}\mathbf{R}{+}\eta\gamma\mathbf{I})\\
&&{+}({\Lambda{+}\Psi})\left(S({}^{b}\mathbf{e}_{q}){-}S(\mathbf{R}^{T}\mathbf{q}_{d})S(\mathbf{e}_{3})\right)\big\}
\end{IEEEeqnarray*}

Consequently by substituting the actual points of the equilibia solutions in $\mathbf{A}$ their eigen-structure can now be studied.

Due to the peculiarity of the control task i.e., we demand pointing stabilization and simultaneously the regulation of the angular velocity about the pointing direction to a desired value, we parametrized the attitude configuration as a spherical pendulum through the unit vector $\mathbf{q}$, because the attitude of the rigid body about the pointing direction is irrelevant.
Moreover, because a spherical pendulum has two rotational degrees of freedom, it holds that $\mathbf{q}^{T}\boldsymbol{\xi}=0$ and $\mathbf{q}^{T}\mathbf{R}{}^{b}\boldsymbol{\omega}=0$, meaning that a spherical pendulum does not have a third component of angular velocity.
In contrast to this, in our system the rotating components of the out-runner motor/propeller assembly do not have the symmetry of the spherical pendulum, and the system has three rotational degrees of freedom.
Therefore, for our system, $\mathbf{q}^{T}\mathbf{R}{}^{b}\boldsymbol{\omega}\neq0$ and the constraint:
\begin{IEEEeqnarray}{rCl}
\mathbf{C}\mathbf{x}=
\begin{bmatrix}
\mathbf{q}^{T}&\mathbf{0}_{1\times3}
\end{bmatrix}
\begin{bmatrix}
\boldsymbol{\xi}\\
\delta^{b}{\boldsymbol{\omega}}
\end{bmatrix}=0\label{eq:constr}
\end{IEEEeqnarray}
in regards to $\boldsymbol{\xi}$ should be satisfied at all times.
Thusly the state vector $\mathbf{x}$ should lie in the null space of $\mathbf{C}\in\mathbb{R}^{1\times6}$.
However if $\mathbf{x}(0)$ satisfies (\ref{eq:constr}) then $\mathbf{x}(t)$ complies with (\ref{eq:constr}), for all $t$, due to the structure of (\ref{eq:cld}) and (\ref{eq:mat}), i.e., (\ref{eq:constr}) and its derivative were embedded in (\ref{eq:mat}) during the linearization procedure.

\section{Eigen-structure of $(\mathbf{q}_{e},{}^{b}\boldsymbol{\omega}_{e})$\label{sec:eig}}

The closed loop properties of (\ref{eq:u}) are investigated by selecting the desired signals to resemble actual reference commands that arise during motor operation on the aerial platform.
Consequently, the resulting reference command is given by,
\begin{IEEEeqnarray}{rCl}
\mathbf{q}_{d}{=}\mathbf{R}_{d}\mathbf{e}_{3},\mathbf{R}_{d}{=}\mathbf{I},{}^{b}\boldsymbol{\omega}_{d}{=}\omega_{d}\mathbf{e}_{3},\omega_{d}{=}1000\frac{rad}{s},{}^{b}\dot{\boldsymbol{\omega}}_{d}{=}\mathbf{0}
\end{IEEEeqnarray}
Using (\ref{eq:cle}), the aforementioned choice corresponds to two equilibrium solutions.
The desired equilibrium,
\begin{IEEEeqnarray}{rCl}
(\mathbf{q}_{e},{}^{b}\boldsymbol{\omega}_{e})=(\mathbf{e}_{3},\omega_{d}\mathbf{e}_{3}),\omega_{d}=1000[{rad}/{s}]\label{eq:des}
\end{IEEEeqnarray}
and the antipodal equilibrium,
\begin{IEEEeqnarray}{rCl}
(\mathbf{q}_{e},{}^{b}\boldsymbol{\omega}_{e})=(-\mathbf{e}_{3},-\omega_{d}\mathbf{e}_{3}),\omega_{d}=1000[{rad}/{s}]\label{eq:anti}
\end{IEEEeqnarray}
The gains of the controller are chosen as in (\ref{eq:gains}).

\subsection{Desired Equilibrium\label{sec:des}}

Using (\ref{eq:des}) with (\ref{eq:mat}) the eigenvalues $\lambda_{i}$, and their corresponding eigenvectors $\mathbf{v}_{i}$, are calculated using MATLAB as,
\begin{IEEEeqnarray}{C}
\lambda_{1}=(-0.0141 + 1.0055i){\cdot}10^{3},\IEEEnonumber\\
\mathbf{v}_{1}=0.0007i\mathbf{e}_{1}+0.0007\mathbf{e}_{2}-0.7071\mathbf{e}_{4}+0.7071i\mathbf{e}_{5},\IEEEnonumber\\
\lambda_{2}=(-2.5693 + 0.0055i){\cdot}10^{3},\IEEEnonumber\\
\mathbf{v}_{2}=-0.0003i\mathbf{e}_{1}+0.0003\mathbf{e}_{2}+0.7071i\mathbf{e}_{4}-0.7071\mathbf{e}_{5},\IEEEnonumber\\
\lambda_{3}=0,\mathbf{v}_{3}=\mathbf{e}_{3},\lambda_{4}=-0.5{\cdot}10^{3},\mathbf{v}_{4}=\mathbf{e}_{6}\IEEEnonumber\\
\lambda_{5}=\bar{\lambda}_{1},\mathbf{v}_{5}=\bar{\mathbf{v}}_{1},\lambda_{6}=\bar{\lambda}_{2},\mathbf{v}_{6}=\bar{\mathbf{v}}_{2}
\end{IEEEeqnarray}
where $\mathbf{e}_{i}\in\mathbb{R}^{6}$, is a base element of the Euclidean space and we have two complex conjugate pairs of eigenvalues $(\lambda_{1,5},\lambda_{2,6})$ and two real eigenvalues $(\lambda_{3},\lambda_{4})$ resulting to six linearly independent associated eigenvectors (two complex conjugate pairs $(\mathbf{v}_{1,5},\mathbf{v}_{2,6})$ and two real $(\mathbf{v}_{3},\mathbf{v}_{4})$).

The base of the null space of (\ref{eq:constr}) is given by,
\begin{IEEEeqnarray}{C}
\mathcal{N}(\mathbf{C})=\{\mathbf{e}_{1},\mathbf{e}_{2},\mathbf{e}_{4},\mathbf{e}_{5},\mathbf{e}_{6}\}\label{eq:null}
\end{IEEEeqnarray}

For $a_{k}\in\mathbb{R},c_{k}\in\mathbb{C}$ the solution of (\ref{eq:mat}) is given by, \cite{Arnold},
\begin{IEEEeqnarray}{C}
\mathbf{x}(t)=\sum_{k=3}^{4}a_{k}e^{\lambda_{k}t}\mathbf{v}_{k}+\sum_{k=1}^{2}c_{k}e^{\lambda_{k}t}\mathbf{v}_{k}+\bar{c}_{k}e^{\bar{\lambda}_{k}t}\bar{\mathbf{v}}_{k}\label{eq:sol}
\end{IEEEeqnarray}
Inspecting (\ref{eq:null}) it is clear that the eigenvector associated with $\lambda_{3}$ does not satisfy (\ref{eq:constr}) since it does not belong to the linear span of (\ref{eq:null}).
Resultantly $a_{3}=0$, $\forall$ $\mathbf{x}(0)$ that satisfy (\ref{eq:constr}).
Thus $\lambda_{3}$ does not partake in (\ref{eq:sol}).
The study of the behavior of $\mathbf{x}(t)$ is not difficult since we have $\text{Re}[\lambda_{2},\lambda_{6}]<\lambda_{4}<\text{Re}[\lambda_{1},\lambda_{5}]<0$.
Therefore the equilibrium (\ref{eq:des}) is an asymptotically stable focus.
Namely we have a rotation with a faster contraction in the ($\mathbf{v}_{2},\mathbf{v}_{6}$)-plane, a fast contraction in the direction of $\mathbf{v}_{4}$ and a rotation with a slower contraction in the ($\mathbf{v}_{1},\mathbf{v}_{5}$)-plane.

\subsection{Antipodal Equilibrium\label{sec:anti}}

Using (\ref{eq:anti}) with (\ref{eq:mat}), the eigenvalues $\lambda_{i}$, and their corresponding eigenvectors $\mathbf{v}_{i}$, are calculated using MATLAB as,
\begin{IEEEeqnarray}{rCl}
\lambda_{1}&=&(-0.7836 + 0.7513i){\cdot}10^{3},\lambda_{5}=\bar{\lambda}_{1}\IEEEnonumber\\
\mathbf{v}_{1}&=&(-0.0005 + 0.0005i)\mathbf{e}_{1}+(0.0005 + 0.0005i)\mathbf{e}_{2}\IEEEnonumber\\
&&- 0.7071i\mathbf{e}_{4}+0.7071\mathbf{e}_{5},\mathbf{v}_{5}=\bar{\mathbf{v}}_{1}\IEEEnonumber\\
\lambda_{2}&=&(2.3670 + 0.2487i){\cdot}10^{3},\lambda_{6}=\bar{\lambda}_{2}\IEEEnonumber\\
\mathbf{v}_{2}&=&0.0003\mathbf{e}_{1}{-}0.0003i\mathbf{e}_{2}{+}0.7071\mathbf{e}_{4}{+}0.7071i\mathbf{e}_{5},\mathbf{v}_{6}=\bar{\mathbf{v}}_{2}\IEEEnonumber\\
\lambda_{3}&=&0,\mathbf{v}_{3}=\mathbf{e}_{3},\lambda_{4}=-0.5{\cdot}10^{3},\mathbf{v}_{4}=\mathbf{e}_{6}
\end{IEEEeqnarray}
We have two complex conjugate pairs of eigenvalues $(\lambda_{1,5},\lambda_{2,6})$ and two real eigenvalues $(\lambda_{3},\lambda_{4})$ resulting to six linearly independent associated eigenvectors (two complex conjugate pairs $(\mathbf{v}_{1,5},\mathbf{v}_{2,6})$ and two real $(\mathbf{v}_{3},\mathbf{v}_{4})$).

The base of the null space of (\ref{eq:constr}) is again given by (\ref{eq:null}), and for $a_{k}\in\mathbb{R},c_{k}\in\mathbb{C}$ the solution of (\ref{eq:mat}) is given by, (\ref{eq:sol}).
Inspecting (\ref{eq:null}) it is clear that the eigenvector associated with $\lambda_{3}$ does not satisfy (\ref{eq:constr}) since it does not belong to the linear span of (\ref{eq:null}).
Resultantly $a_{3}=0$, $\forall$ $\mathbf{x}(0)$ that satisfy (\ref{eq:constr}).
Again $\lambda_{3}$ does not partake in (\ref{eq:sol}).

To study the behavior of $\mathbf{x}(t)$ we note that since we have $\text{Re}[\lambda_{1},\lambda_{5}]{<}\lambda_{4}{<}0{<}\text{Re}[\lambda_{2},\lambda_{6}]$, the equilibrium (\ref{eq:des}) is a saddle point.
Specifically we have a rotation with a faster contraction in the ($\mathbf{v}_{1},\mathbf{v}_{5}$)-plane, a fast contraction in the direction of $\mathbf{v}_{4}$ and a rotation with a fast dilation in the ($\mathbf{v}_{2},\mathbf{v}_{6}$)-plane.

\section{Precession/Nutation Frequency Estimation}

In Section \ref{sec:pdav}, studying a typical PDAV trajectory, high frequency nutation/precession oscillations were observed.
The emergence of these non trivial nutation/precession dynamic oscillations (see Fig. \ref{fig:pdav_3} (top and middle rows)), are of significant importance, since suitable pointing actuators need to be able to handle/negotiate such high frequency oscillations.
Thus a systematic method of estimating the frequency of the nutation/precession oscillations is needed.

Furthermore, the ability to obtain estimates of the nutation/precession oscillations frequency provides a tool in gaining understanding on the PDAV closed-loop intricacies, a tool for actuator selection and finally a criterion on the feasibility of an experimental implementation, i.e., if for a required thrust no suitable pointing actuators exist or the actuators are bulky for the needed application, then other means of thrust generation should be considered.

Our approach in developing this tool originates from analyzing the process of trajectory tracking.
During trajectory tracking, initially the states of the system are at equilibrium and by gradually/smoothly shifting the reference PDAV command, the global equilibrium of the system changes to a new preferred state.
This procedure can be thought off roughly as a concatenation of infinitesimal flows (solutions) chasing an infinitesimally shifting reference PDAV command.
Thus, we assume that a valid estimate of the frequency of the nutation/precession oscillations can be obtained, by using the solutions of the linearized system, (\ref{eq:sol}), in the neighborhood of the desired equilibrium.
Furthermore this frequency estimate can be extended to the entire duration of the trajectory tracking maneuver, if the maneuver is performed in a sufficiently smooth and gradual manner.
This assumption guides our analysis with the developed procedure to follow next.

The solution of (\ref{eq:mat}), for $a_{k}\in\mathbb{R},c_{k}\in\mathbb{C}$, as given by (\ref{eq:sol}) can be rewritten in the form, \cite{Arnold},
\begin{IEEEeqnarray}{C}
\mathbf{x}(t)=\sum_{k=3}^{4}a_{k}e^{\lambda_{k}t}\mathbf{v}_{k}+2\text{Re}\left[\sum_{k=1}^{2}c_{k}e^{\lambda_{k}t}\mathbf{v}_{k}\right]\label{eq:solre}
\end{IEEEeqnarray}

The eigenvalues, $\lambda_{i}$, and eigenvectors, $\mathbf{v}_{i}$, from Section \ref{sec:des} and the fact that $a_{3}=0$ $\forall$ $\mathbf{x}(0)$ that satisfy (\ref{eq:constr}), are applied to (\ref{eq:solre}) to find the solutions of the linearized system in the neighborhood of the equilibrium.
The result is,
\begin{IEEEeqnarray}{rCl}
\mathbf{x}(t)&=&a_{4}e^{\lambda_{4}t}\mathbf{e}_{6}+2\sum_{k=1}^{2}\Big\{\text{Re}[c_{k}](\cos(\mu_{k}t)\text{Re}[\mathbf{v}_{k}]\IEEEnonumber\\
&&-\sin(\mu_{k}t)\text{Imag}[\mathbf{v}_{k}])-\text{Imag}[c_{k}](\cos(\mu_{k}t)\text{Imag}[\mathbf{v}_{k}]\IEEEnonumber\\
&&+\sin(\mu_{k}t)\text{Re}[\mathbf{v}_{k}])\Big\}\label{eq:solreexp}
\end{IEEEeqnarray}
where $\pi_{k},\mu_{k}$ are obtained from,
\begin{IEEEeqnarray}{rCl}
\lambda_{k}=\pi_{k}+\mu_{k}i\IEEEnonumber
\end{IEEEeqnarray}
Using the ''313'' Euler sequence, the expression that correlates the precession rate, $\dot{\phi}$ and nutation rate, $\dot{\theta}$ with the angular velocity ${}^{b}\boldsymbol{\omega}$ of the system is given by,
\begin{IEEEeqnarray}{rCl}
{\begin{bmatrix}
\dot{\phi}\\
\dot{\theta}\\
\dot{\psi}
\end{bmatrix}}{=}\begin{bmatrix}
\frac{\sin(\psi)}{\sin(\theta)}&\frac{\cos(\psi)}{\sin(\theta)}&0\\
\cos(\psi)&-\sin(\psi)&0\\
-\frac{\sin(\psi)\cos(\theta)}{\sin(\theta)}&-\frac{\cos(\psi)\cos(\theta)}{\sin(\theta)}&1
\end{bmatrix}{\begin{bmatrix}
{}^{b}\boldsymbol{\omega}_{1}\\
{}^{b}\boldsymbol{\omega}_{2}\\
{}^{b}\boldsymbol{\omega}_{3}
\end{bmatrix}}
\label{eq:omegaeuler}
\end{IEEEeqnarray}
where $\psi$ is the third Euler angle denoting the rotation angle of the out-runner shell/propeller about its own axis $\mathbf{e}_{3}$.
Note that the expression of the precession rate $\dot{\phi}$, given by (\ref{eq:omegaeuler}) is only valid when $\theta{\notin} \{0,k\pi\lvert k{\in}\mathbb{Z}\}$ i.e., the propeller axis, $\mathbf{e}_{3}$, is not parallel with the vertical, $\mathbf{E}_{3}$.
This agrees with the physical intuition that no precession exists when the motor is aligned with $\mathbf{E}_{3}$.
The sixth component of (\ref{eq:solreexp}) equals to,
\begin{IEEEeqnarray}{rCl}\mathbf{x}(t)\cdot\mathbf{e}_{6}&=&\delta^{b}{\boldsymbol{\omega}}_{3}=a_{4}e^{\lambda_{4}t}=a_{4}e^{-0.5{\cdot}10^{3}t}\label{eq:linspin}
\end{IEEEeqnarray}

Equation (\ref{eq:linspin}) with (\ref{eq:delo}) imply that in the neighborhood of the desired equilibrium ${}^{b}{\boldsymbol{\omega}}_{3}{\approx}\omega_{d}$ because the time constant obtained from (\ref{eq:linspin}) is extremely small ($\tau=0.002$s).
Thus the propeller angular speed is regulated to the desired value very fast.
Since ${}^{b}\boldsymbol{\omega}_{3}=\dot{\psi}+\dot{\phi}\cos(\theta)$ and $\dot{\psi}\ggg\dot{\phi}$, it holds that ${}^{b}\boldsymbol{\omega}_{3}\approx\dot{\psi}=\omega_{d}$.
Resultantly the rotation angle of the propeller $\psi(t)$ in the neighborhood of the desired equilibrium is,
\begin{IEEEeqnarray}{rCL}
\psi(t)\approx{}^{b}\boldsymbol{\omega}_{3}t+\psi(t=0)=\omega_{d}t+\psi_{0}
\label{eq:psiaprox}
\end{IEEEeqnarray}
Substituting (\ref{eq:psiaprox}) and the components $\delta^{b}{\boldsymbol{\omega}}_{1}=\mathbf{x}(t)\cdot\mathbf{e}_{4}$, $\delta^{b}{\boldsymbol{\omega}}_{2}=\mathbf{x}(t)\cdot\mathbf{e}_{5}$ in (\ref{eq:omegaeuler}), for the nutation rate $\dot{\theta}$ we get,
\begin{IEEEeqnarray}{rCl}
\dot{\theta}&=&
2\cos(\omega_{d}t+\psi_{0})\sum_{k=1}^{2}\mathbf{e}_{4}\cdot\Big\{\text{Re}[c_{k}](\cos(\mu_{k}t)\text{Re}[\mathbf{v}_{k}]\IEEEnonumber\\
&&-\sin(\mu_{k}t)\text{Imag}[\mathbf{v}_{k}])-\text{Imag}[c_{k}](\cos(\mu_{k}t)\text{Imag}[\mathbf{v}_{k}]\IEEEnonumber\\
&&+\sin(\mu_{k}t)\text{Re}[\mathbf{v}_{k}])\Big\}\IEEEnonumber\\
&&-2\sin(\omega_{d}t+\psi_{0})\sum_{k=1}^{2}\mathbf{e}_{5}\cdot\Big\{\text{Re}[c_{k}](\cos(\mu_{k}t)\text{Re}[\mathbf{v}_{k}]\IEEEnonumber\\
&&-\sin(\mu_{k}t)\text{Imag}[\mathbf{v}_{k}])-\text{Imag}[c_{k}](\cos(\mu_{k}t)\text{Imag}[\mathbf{v}_{k}]\IEEEnonumber\\
&&+\sin(\mu_{k}t)\text{Re}[\mathbf{v}_{k}])\Big\}\label{eq:inprecest}
\end{IEEEeqnarray}

	To simplify (\ref{eq:inprecest}) we inspect  the eigenvalues of Section \ref{sec:des} and observe that $\text{Re}[\lambda_{2}]{=}\pi_{2}{=}{-}2569.3{\lll}\text{Re}[\lambda_{1}]{=}\pi_{1}{=}{-}14.1$.
Resultantly the solution contracts extremely fast in the  $(\mathbf{v}_{2},\mathbf{v}_{6})$-plane and the nutation rate can be approximated as,
\begin{IEEEeqnarray}{rCl}
\dot{\theta}&\approx&
2\cos(\omega_{d}t)\Big\{\breve{C}_{1}\cos(\mu_{1}t)+\breve{D}_{1}\sin(\mu_{1}t)\Big\}\IEEEnonumber\\
&&-2\sin(\omega_{d}t)\Big\{\breve{C}_{2}\cos(\mu_{1}t)+\breve{D}_{2}\sin(\mu_{1}t)\Big\}\\
\breve{C}_{1}&=&\cos(\psi_{0})C_{1}{-}\sin(\psi_{0})C_{2},\breve{D}_{1}{=}\cos(\psi_{0})D_{1}-\sin(\psi_{0})D_{2}\IEEEnonumber\\
\breve{C}_{2}&=&\cos(\psi_{0})C_{2}+\sin(\psi_{0})C_{1},\breve{D}_{2}{=}\cos(\psi_{0})D_{2}+\sin(\psi_{0})D_{1}
\IEEEnonumber\\
C_{i}&=&\text{Re}[c_{1}]\text{Re}[\mathbf{v}_{1_{i+3}}]-\text{Imag}[c_{1}]\text{Imag}[\mathbf{v}_{1_{i+3}}],i=1,2\IEEEnonumber\\
D_{i}&=&-\text{Re}[c_{1}]\text{Imag}[\mathbf{v}_{1_{i+3}}]-\text{Imag}[c_{1}]\text{Re}[\mathbf{v}_{1_{i+3}}],i=1,2\IEEEnonumber
\end{IEEEeqnarray}
Employing product-to-sum identities and rearranging we get,
\begin{IEEEeqnarray}{rCl}
\dot{\theta}&\approx&(\breve{C}_{1}-\breve{D}_{2})\cos((\omega_{d}{+}\mu_{1})t)+(\breve{D}_{1}-\breve{C}_{2})\sin((\omega_{d}{+}\mu_{1})t)\IEEEnonumber\\
&&+(\breve{C}_{1}+\breve{D}_{2})\cos((\omega_{d}{-}\mu_{1})t)-(\breve{D}_{1}+\breve{C}_{2})\sin((\omega_{d}{-}\mu_{1})t)\IEEEnonumber
\end{IEEEeqnarray}
Additionally, using (\ref{eq:mat}) for several values of $\omega_{d}$ we observed that $\mu_{1}\approx\omega_{d}$.
Resultantly the terms $\omega_{d}+\mu_{1}$ dominate the oscillation frequency, $f_{n}$, since $\omega_{d}+\mu_{1}\ggg\omega_{d}-\mu_{1}$ and the final estimation of the nutation oscillation frequency is,
\begin{IEEEeqnarray}{rCl}
f_{n}\approx\frac{\omega_{d}{+}\mu_{1}}{2\pi}\label{eq:freqest}
\end{IEEEeqnarray}

Thus, for the desired command in Section \ref{sec:eig}, the frequency obtained by means of FFT analysis (see Fig. \ref{fig:pdav_3} (bottom row)) gives that $f_{n_{FFT}}=318.6454\text{Hz}$ while using (\ref{eq:freqest}) we obtain $f_{n}=319.18\text{Hz}$ which is almost identical to the measured frequency $f_{n_{FFT}}$.
This validates the developed formula, i.e., (\ref{eq:freqest}), and we are now equipped with a method to estimate the high frequency nutation oscillations during PDAV trajectory tracking in the neighborhood of the equilibrium and thus during a smooth PDAV trajectory.

Note that in the preceding analysis we chose to develop an expression that estimates the high frequency nutation oscillations but the procedure can be repeated using the precession component of (\ref{eq:omegaeuler}) to get an expression similar to (\ref{eq:freqest}).

\section{Flow ''close'' to the PDAV Equilibria}

At the antipodal saddle equilibrium, the conditions of the Hartman-Grobman theorem \cite{Guckenheimer}, and the Stable Manifold theorem \cite{Guckenheimer}, are satisfied.
Thus a \textit{local stable manifold},
\begin{IEEEeqnarray}{l}
\mathbf{W}^{s}_{loc}(-\mathbf{q}_{d},-{}^{b}\boldsymbol{\omega}_{d}){=}\{\mathbf{x}{\in}\mathbf{U}|\lim_{t\to\infty}\boldsymbol{\varphi}^{t}(\mathbf{x}){=}(\ref{eq:anti}),\forall t{\geq}0\}\label{eq:wsloc}
\end{IEEEeqnarray}
exists tangent to the flat stable eigenspace $\mathbf{E}^{s}$ of the linearized system.
Note that $\boldsymbol{\varphi}^{t}$ in (\ref{eq:wsloc}) is the forward flow map.
Furthermore, the global stable manifold, $\mathbf{W}^{s}$, can be obtained by letting points in $\mathbf{W}^{s}_{loc}$ flow backwards in time:

\begin{IEEEeqnarray}{rCl}
\mathbf{W}^{s}(-\mathbf{q}_{d},-{}^{b}\boldsymbol{\omega}_{d})&{=}&\bigcup_{t\geq0}\boldsymbol{\varphi}^{-t}\left(\mathbf{W}^{s}_{loc}(-\mathbf{q}_{d},{}^{b}\boldsymbol{\omega}_{d})\right)\label{eq:stman}
\end{IEEEeqnarray}
where $\mathbf{U}\subset\text{S}^{2}{\times}\mathbb{R}^{3}$, a neighborhood of (\ref{eq:anti}), and $\boldsymbol{\varphi}^{-t}$ is the backward flow map \cite{Guckenheimer}.

The existence of the stable manifold $\mathbf{W}^{s}$ has significant implications on the closed-loop system, since trajectories on it converge to the antipodal equilibrium, while trajectories near it need significant time to converge to the desired equilibrium.

However, we are more interested in the flow ''close'' to the invariant manifold since in an experimental implementation, disturbances during operation ensure that the states will not remain on the invariant manifold.
Additionally during PDAV trajectory tracking, the existence of disturbances can produce a shift of the state close to the invariant manifold.
This event can result in irregular behavior and a goal of this paper is to comprehend the severity of this event.
Thus obtaining an understanding of the flow ''close'' to the stable manifold of $(-\mathbf{q}_{d},-{}^{b}\boldsymbol{\omega}_{d})$ is more useful in our case.

\subsection{Flow ''close'' to the Antipodal Equilibrium\label{sec:fcse}}

We follow the method presented in \cite{Krauskopf}, and later used in \cite{SMSE}:
The stable eigenvectors $\mathbf{v}_{1}$, $\mathbf{v}_{3}$, $\mathbf{v}_{5}$, of (\ref{eq:anti}) from Section \ref{sec:anti} are used to generate the local stable eigenspace $\mathbf{E}^{s}_{loc}$:
\begin{IEEEeqnarray}{l}
\mathbf{E}^{s}_{loc}(-\mathbf{q}_{d},-{}^{b}\boldsymbol{\omega}_{d}){=}\{(\mathbf{q},{}^{b}\boldsymbol{\omega})\in\text{S}^{2}{\times}\mathbb{R}^{3}|\IEEEnonumber\\
\mathbf{q}=\exp\left( S(\Delta_{q}\{\varepsilon\cos(\vartheta)(\sigma\mathbf{v}_{1}+\bar{\sigma}\mathbf{v}_{5})\}\right)({-}\mathbf{q}_{d}),\sigma\in\mathbb{C},\IEEEnonumber\\
{}^{b}\boldsymbol{\omega}=-{}^{b}\boldsymbol{\omega}_{d}+\Delta_{\omega}(\varepsilon\cos(\vartheta)(\sigma\mathbf{v}_{1}{+}\bar{\sigma}\mathbf{v}_{5})+\varsigma\sin(\vartheta)\mathbf{v}_{3}),\IEEEnonumber\\
\varepsilon,\varsigma{\lll}1,\vartheta{\in}[0,2\pi),\Delta_{q}{=}[\mathbf{I},\mathbf{0}],\Delta_{\omega}{=}[\mathbf{0},\mathbf{I}]\in\mathbb{R}^{3\times6})\label{eq:elocs}
\end{IEEEeqnarray}
A distance metric on the tangent bundle is defined as:
\begin{IEEEeqnarray}{l}
\text{d}_{\mathbf{q},\omega}((\mathbf{q}_{1},{}^{b}\boldsymbol{\omega}_{1}),(\mathbf{q}_{2},{}^{b}\boldsymbol{\omega}_{2})){=}\Psi(\mathbf{q}_{1},\mathbf{q}_{2}){+}\lVert{}^{b}\boldsymbol{\omega}_{1}{-}{}^{b}\boldsymbol{\omega}_{2}\rVert\label{eq:dmet}
\end{IEEEeqnarray}
This distance metric will be used to check if the method proposed in \cite{SMSE} for visualizing the stable manifold can also be used for (\ref{eq:cld}) and additionally as a measure of proximity to (\ref{eq:anti}).
The backward flow map $\boldsymbol{\varphi}^{-t}$ and forward flow map $\boldsymbol{\varphi}^{t}$ are calculated using variational integrators \cite{SMSE}, \cite{Orbit}.

\begin{figure}[!h]
\centering
\subfloat[$t=0.0366$.\label{fig:saddle_1}]{\includegraphics[width=0.49\columnwidth]{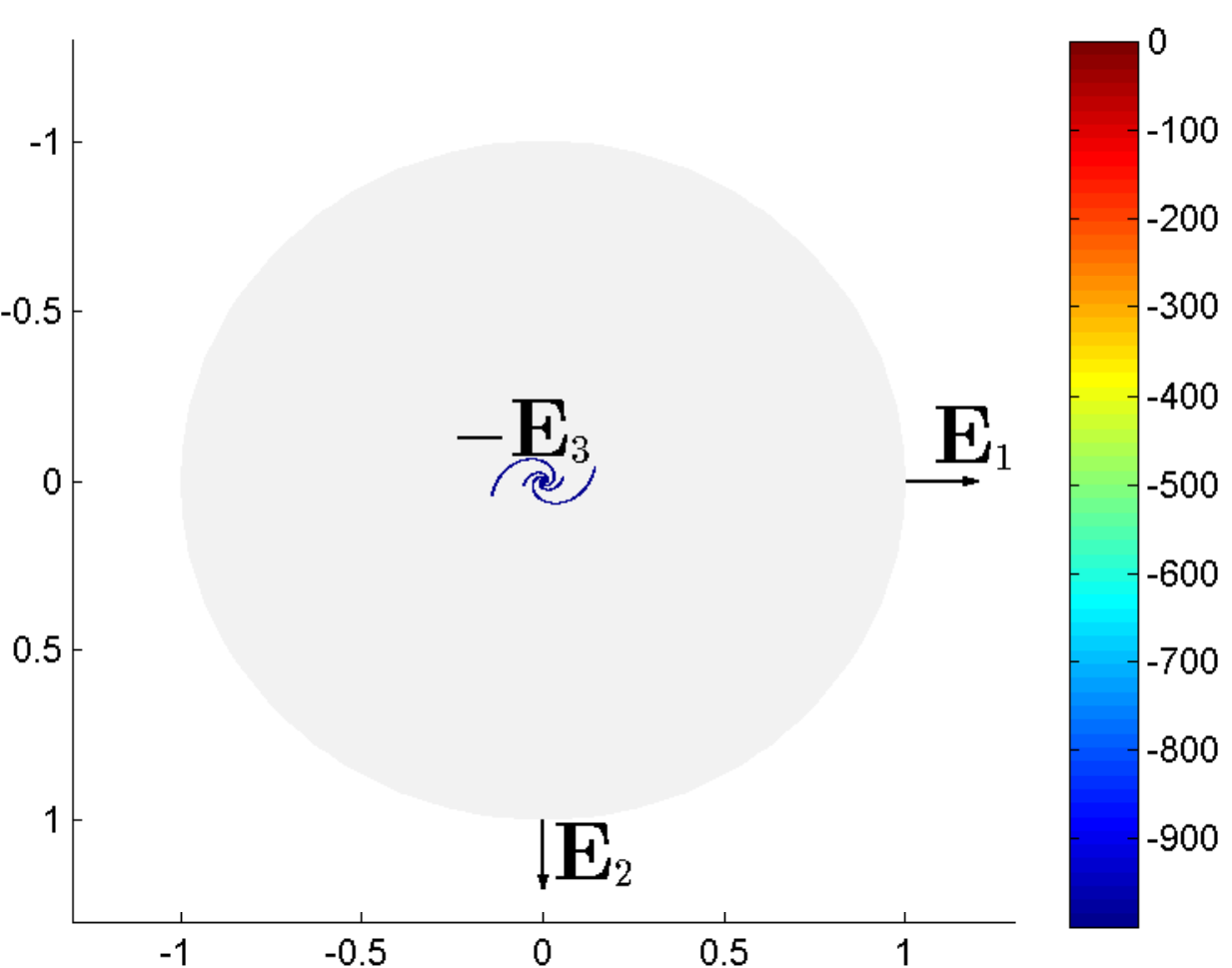}}
~
\subfloat[$t=0.0409$.\label{fig:saddle_2}]{\includegraphics[width=0.49\columnwidth]{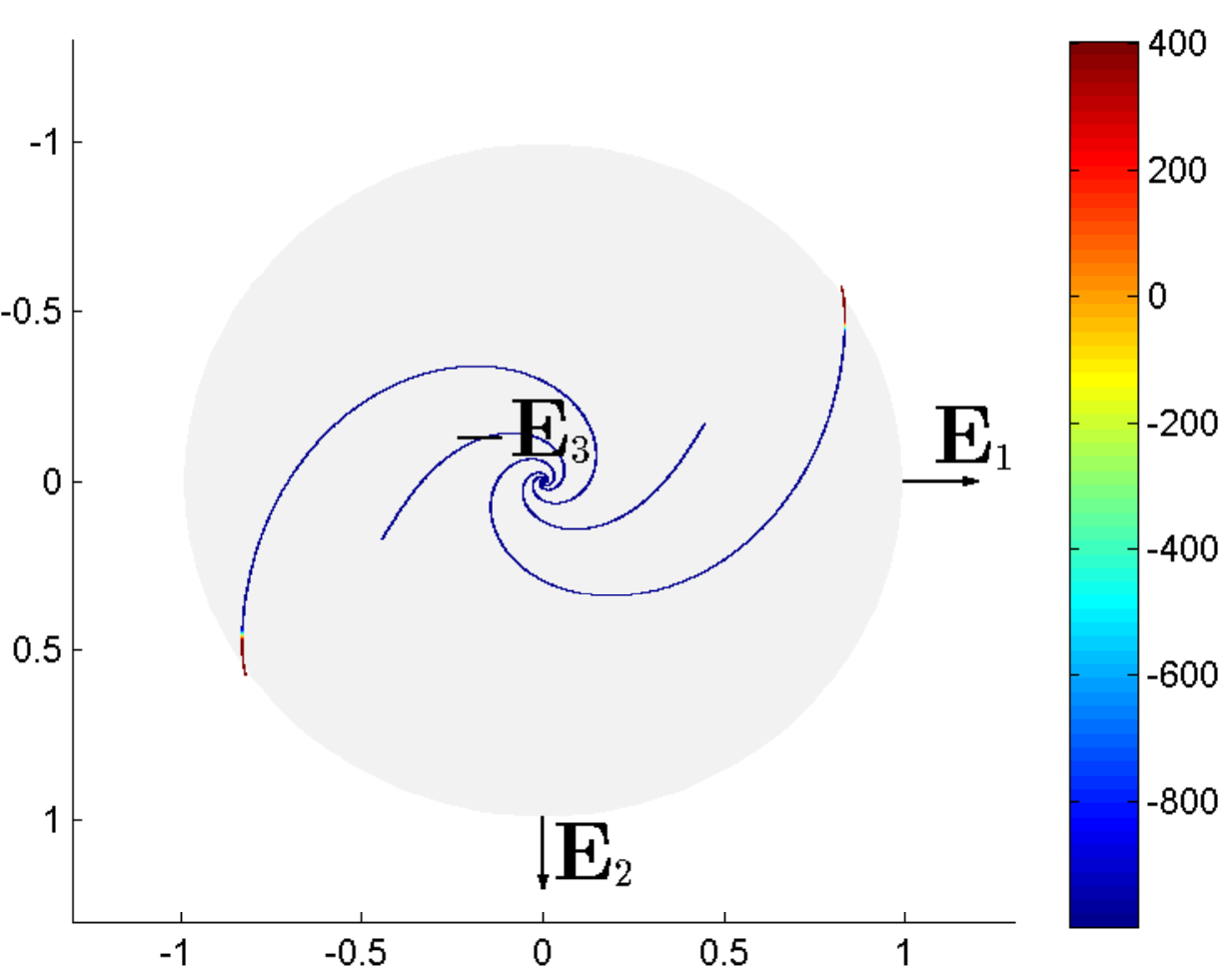}}

\subfloat[$t=0.0420$.\label{fig:saddle_3}]{\includegraphics[width=0.49\columnwidth]{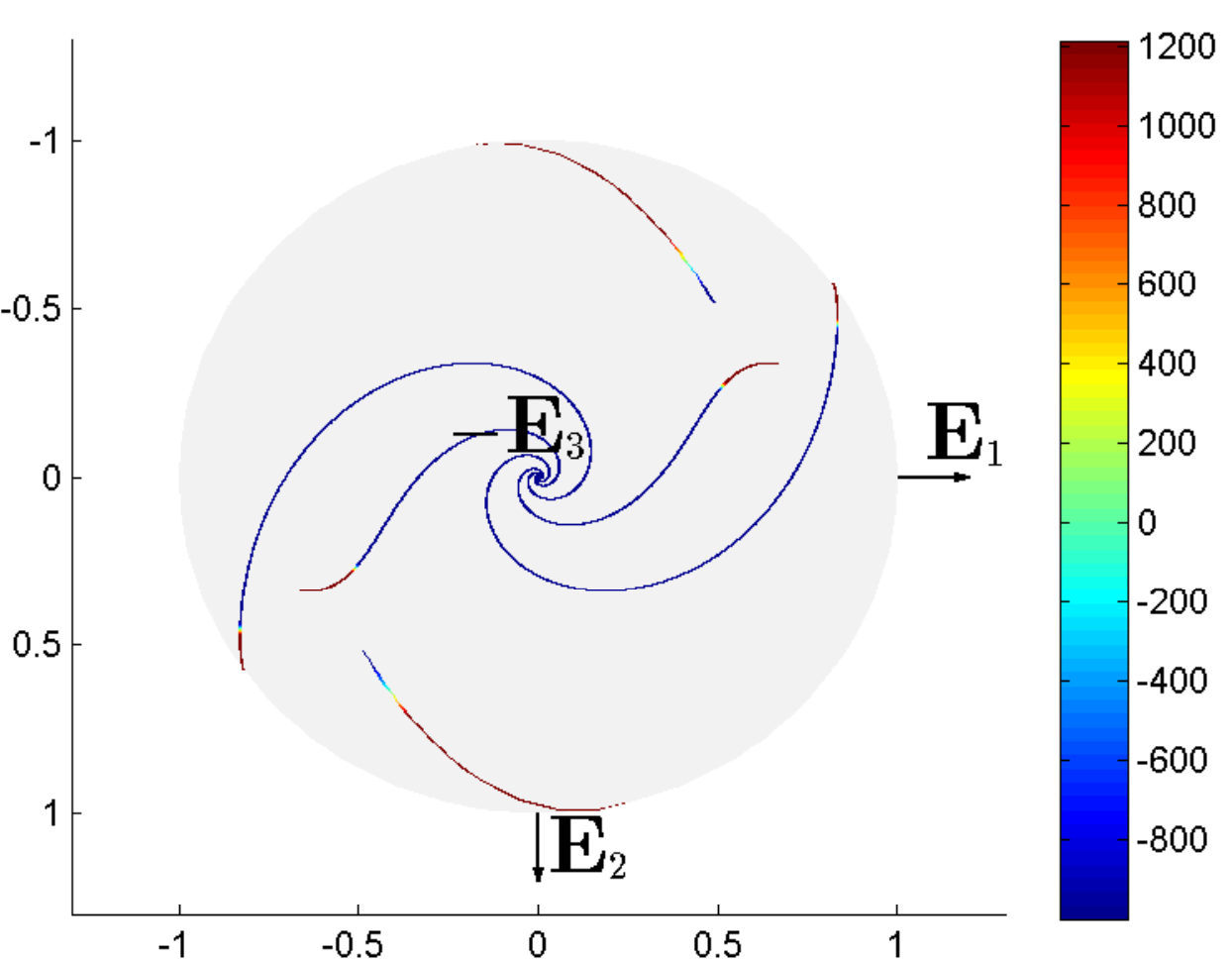}}
~
\subfloat[$t=0.0433$.\label{fig:saddle_4}]{\includegraphics[width=0.49\columnwidth]{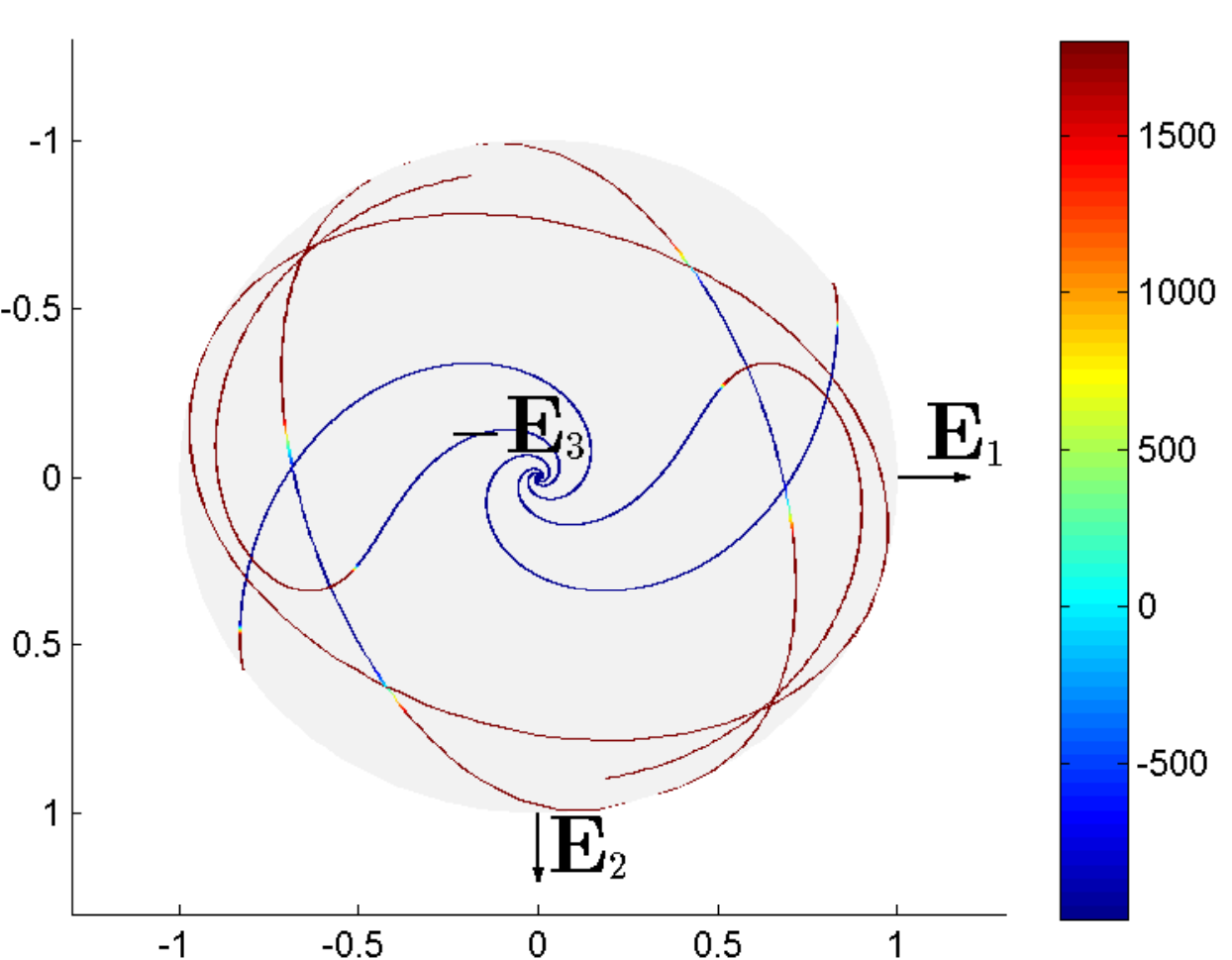}}

\subfloat[$t=0.0450$.\label{fig:saddle_5}]{\includegraphics[width=0.49\columnwidth]{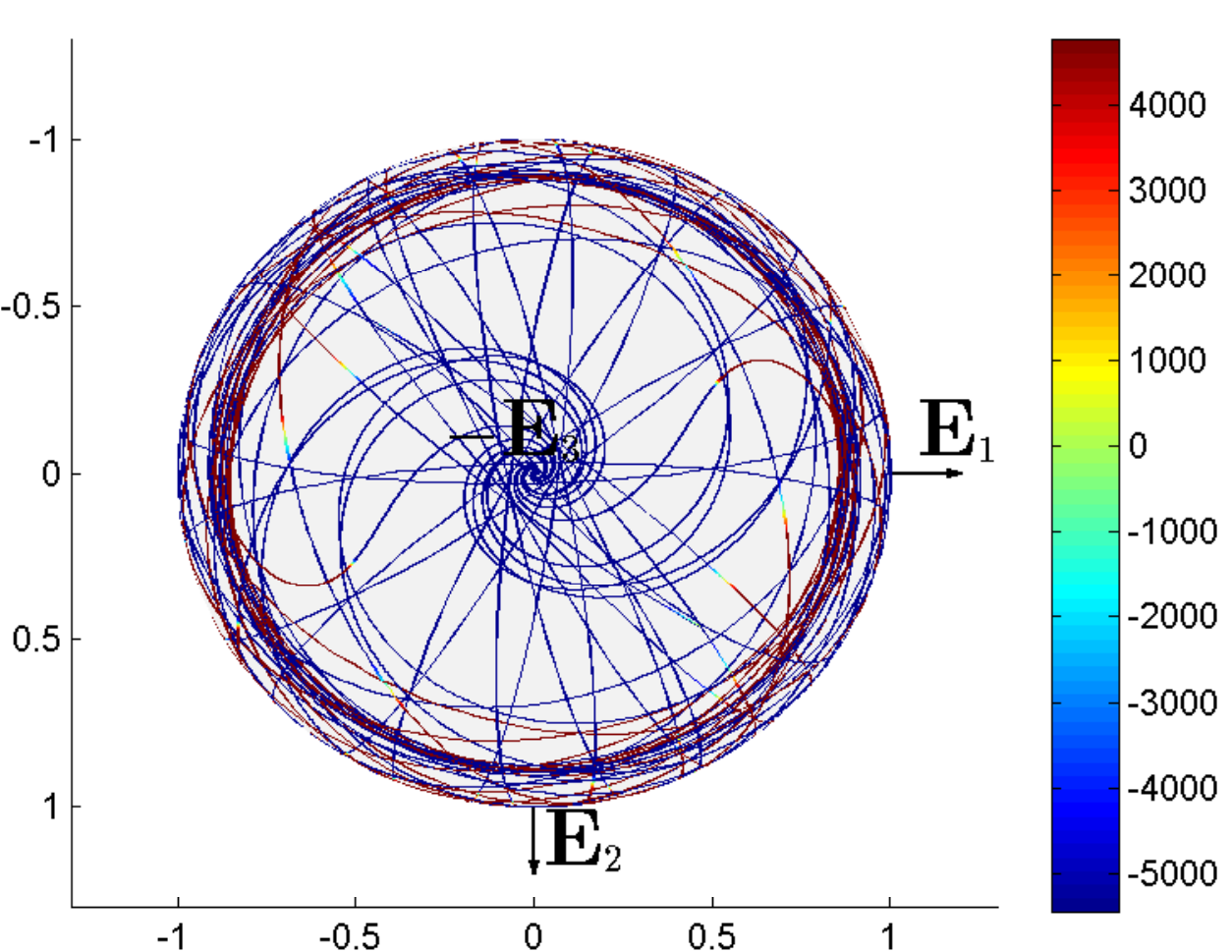}}
~
\subfloat[$t=0.0450$.\label{fig:saddle_6}]{\includegraphics[width=0.49\columnwidth]{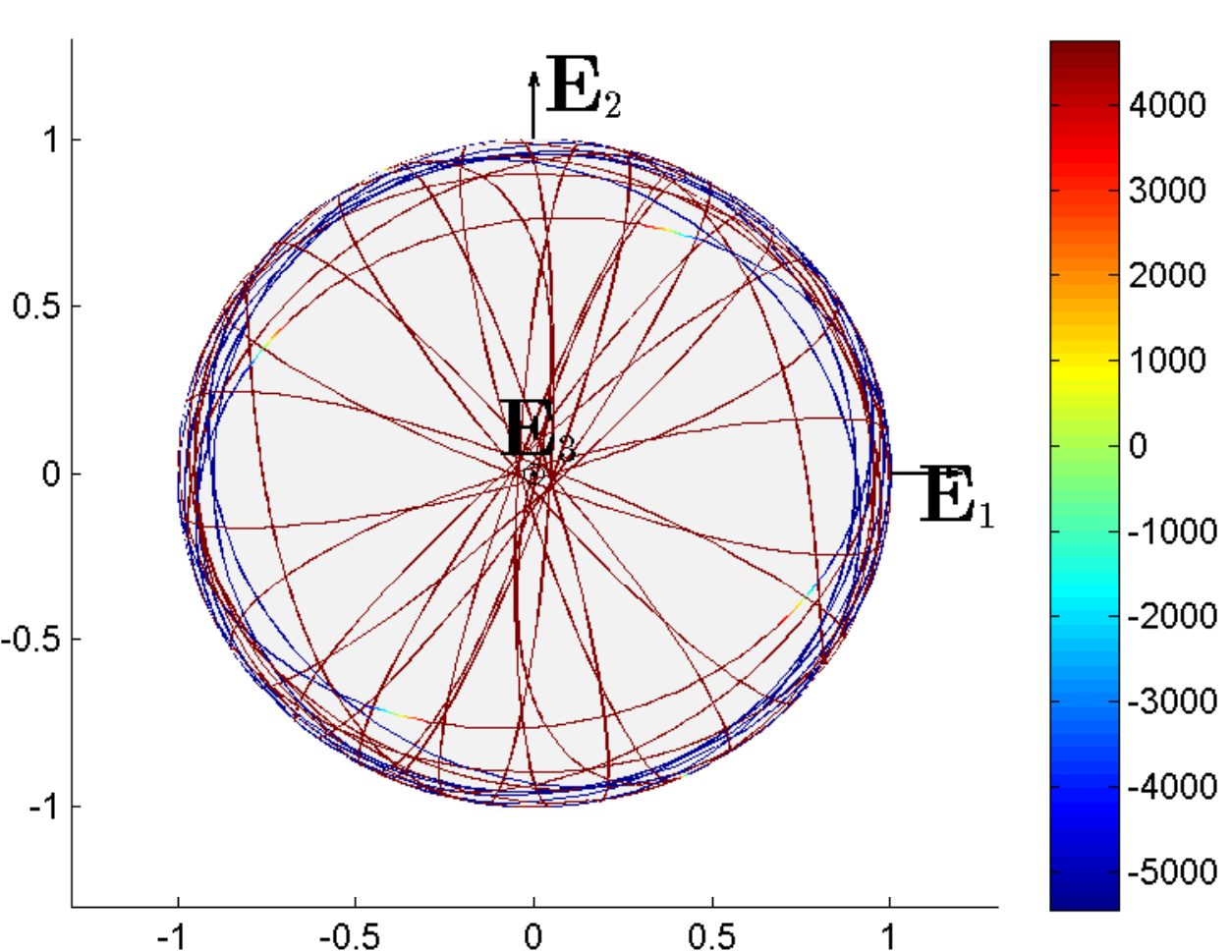}}

\caption{Backwards flow represented by $\{\boldsymbol{\varphi}^{{-}t}(\mathbf{p})\}_{t{>}0}$, $\mathbf{p}{\in}(\ref{eq:elocs})$ extremely ''close'' to the stable manifold of the saddle $({-}\mathbf{q}_{d},-{}{}^{b}\boldsymbol{\omega}_{d})$.
The trajectories are depicted on $\text{S}^{2}$ and the magnitude of the angular velocity about the pointing direction is indicated by the color of the trajectories according to the colorbar.
(\ref{fig:saddle_1}-\ref{fig:saddle_4}) Four points from (\ref{eq:elocs}), integrated backwards in time and shown for instances of $t$.
(\ref{fig:saddle_5}-\ref{fig:saddle_6}) Flow shown in both sides of $\text{S}^{2}$ for ten points.
}
\label{fig:saddle}
\end{figure}

We pick ten points from (\ref{eq:elocs}), with $\varepsilon{=}\varsigma{=}1{\cdot}10^{-6}$, $\sigma{=}1{+}1i$.
The trajectories evolving on $\text{S}^{2}$ are shown in Fig. \ref{fig:saddle} where each colored path stems from one of the ten selected points.
The angular velocity about the pointing direction is indicated by the color of the trajectories according to the colorbar.

Several observations regarding trajectories that converge in  (\ref{eq:elocs}) and come ''close'' to $\mathbf{W}^{s}$ are summarized next:

The trajectories near the saddle are approximately logarithmic spirals that as they move away from (\ref{eq:anti}) they wrap around $\text{S}^{2}$ in an intricate manner.
The angular velocity about the pointing direction converges to $-{}{}^{b}\boldsymbol{\omega}_{d}$ as $\mathbf{q}$ approaches the saddle, even from extreme positive or negative initial spinning velocities (see Fig. \ref{fig:saddle_5}).
Some of the trajectories that come ''close'' to $\mathbf{W}^{s}$ wrap around $\text{S}^{2}$ multiple times.
Resultantly even if $\mathbf{W}^{s}$ is of zero measure, during operation $\mathbf{q}$ might approach arbitrary close to the saddle $(-\mathbf{q}_{d},-{}{}^{b}\boldsymbol{\omega}_{d})$.

\begin{figure}[!h]
\centering
\subfloat[\label{fig:saddle_f_1}]{\includegraphics[width=0.49\columnwidth]{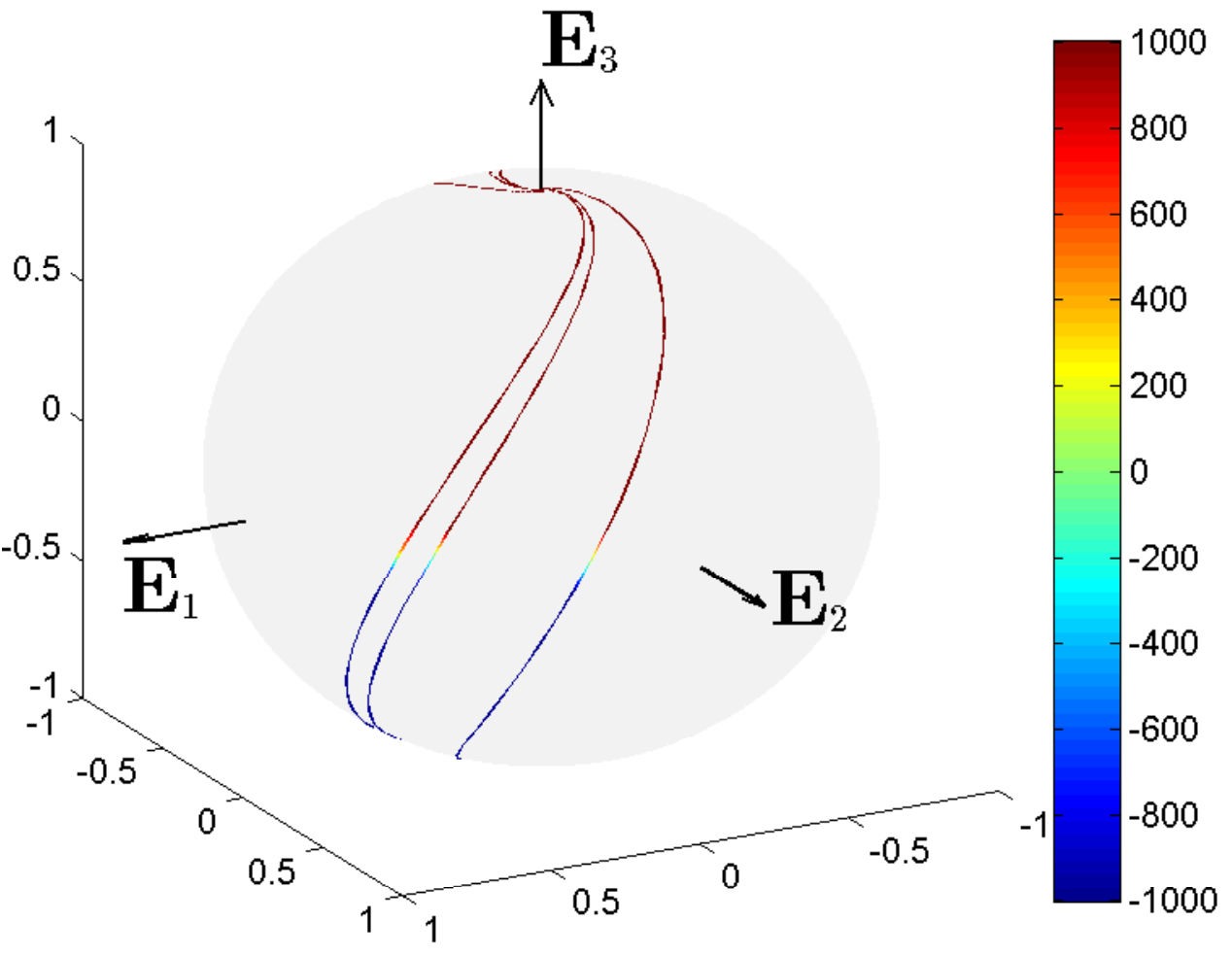}}
~
\subfloat[\label{fig:saddle_f_2}]{\includegraphics[width=0.49\columnwidth]{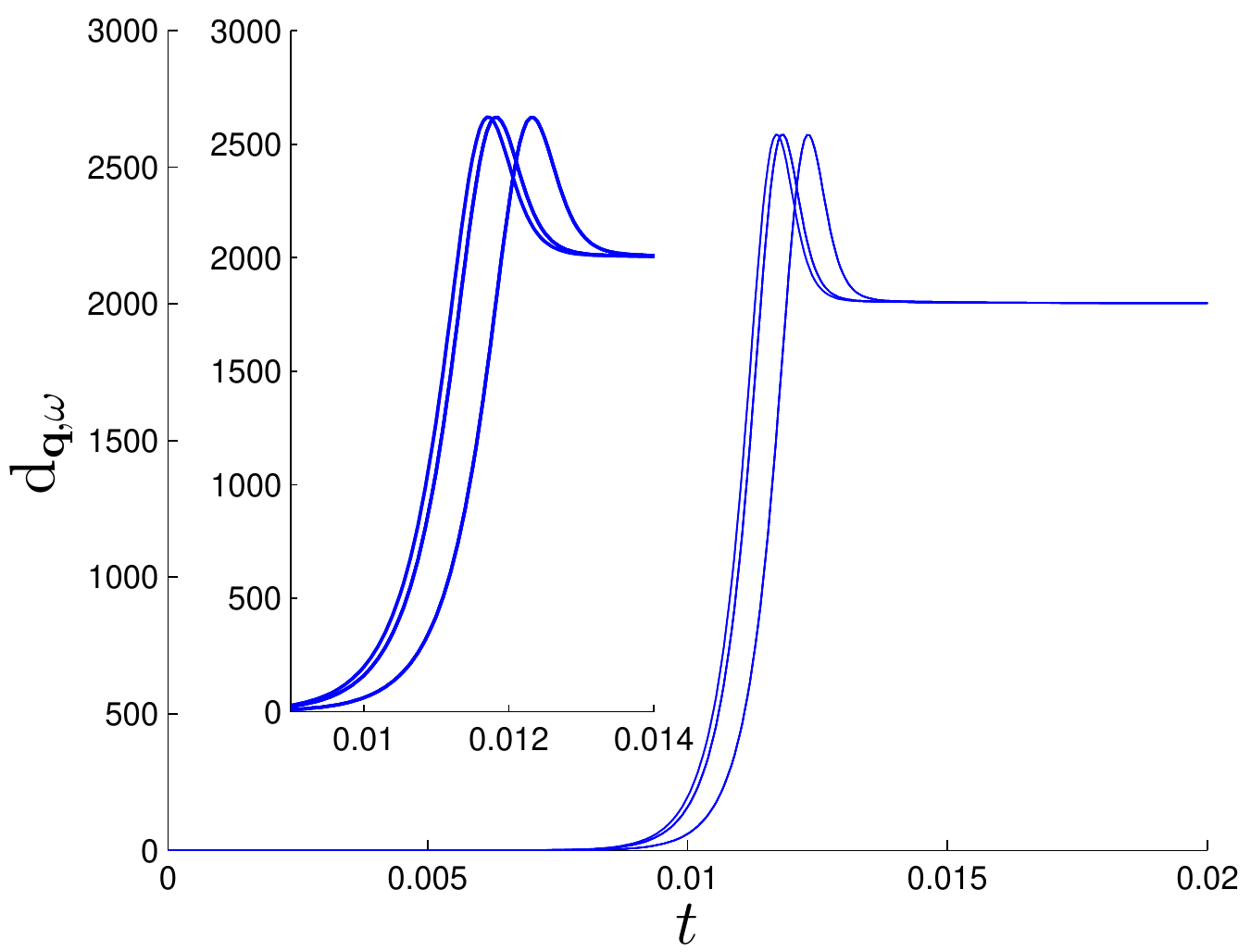}}

\caption{Forward flow represented by $\{\boldsymbol{\varphi}^{t}(\mathbf{p}_{s})\}_{t{>}0}$, for the ten points depicted in Fig. (\ref{fig:saddle}).
(\ref{fig:saddle_f_1}) Trajectories evolving on $\text{S}^{2}$ with ${}^{b}\boldsymbol{\omega}\cdot\mathbf{e}_{3}$ illustrated by the color of the trajectories according to the colorbar.
(\ref{fig:saddle_f_2}) Distance metric $\text{d}_{\mathbf{q},\omega}((-\mathbf{q}_{d},-{}^{b}\boldsymbol{\omega}_{d}),\boldsymbol{\varphi}^{t}(\mathbf{p}_{s}))$ during $\{\boldsymbol{\varphi}^{t}(\mathbf{p}_{s})\}_{t{>}0}$.
}
\label{fig:saddle_f}
\end{figure}

In an attempt to produce $\mathbf{W}^{s}(-\mathbf{q}_{d},-{}{}^{b}\boldsymbol{\omega}_{d})$, several simulations were conducted using $\varepsilon\lll1{\cdot}10^{-6}$, $\varsigma\lll1{\cdot}10^{-6}$ aiming to find points, $\mathbf{p}_{s}$, in $\mathbf{W}^{s}_{loc}(-\mathbf{q}_{d},-{}^{b}\boldsymbol{\omega}_{d})$.
To check that the produced points actually belong to $\mathbf{W}^{s}_{loc}$, the forward flow map $\boldsymbol{\varphi}^{t}$ and (\ref{eq:dmet}) were used to check if $\boldsymbol{\varphi}^{t}(\mathbf{p}_{s})\to (-\mathbf{q}_{d},-{}^{b}\boldsymbol{\omega}_{d})$ as $t\to\infty$.
The flow is shown in Fig. \ref{fig:saddle_f_1}, and the values of $\text{d}_{\mathbf{q},\omega}((-\mathbf{q}_{d},-{}^{b}\boldsymbol{\omega}_{d}),\boldsymbol{\varphi}^{t}(\mathbf{p}_{s}))$ in Fig. \ref{fig:saddle_f_2}, both indicating that despite approaching very ''close'' to $(-\mathbf{q}_{d},-{}^{b}\boldsymbol{\omega}_{d})$, the flow eventually converges to $(\mathbf{q}_{d},{}^{b}\boldsymbol{\omega}_{d})$.
Thus despite the fact that $\mathbf{E}^{s}_{loc}(-\mathbf{q}_{d},-{}^{b}\boldsymbol{\omega}_{d})$ is tangent to $\mathbf{W}^{s}_{loc}(-\mathbf{q}_{d},-{}^{b}\boldsymbol{\omega}_{d})$, points calculated using $\mathbf{E}^{s}_{loc}$ are ''close'' but do not necessarily belong to $\mathbf{W}^{s}_{loc}$.
Therefore, the method proposed in \cite{SMSE} for the calculation of $\mathbf{W}^{s}(-\mathbf{q}_{d},-{}{}^{b}\boldsymbol{\omega}_{d})$ by $\boldsymbol{\varphi}^{-t}$ does not work here for (\ref{eq:cld}).

\subsection{Flow ''close'' to the Desired Equilibrium}

\begin{figure}[!h]
\centering
\subfloat[$t=0.0058$.\label{fig:focus_1}]{\includegraphics[width=0.49\columnwidth]{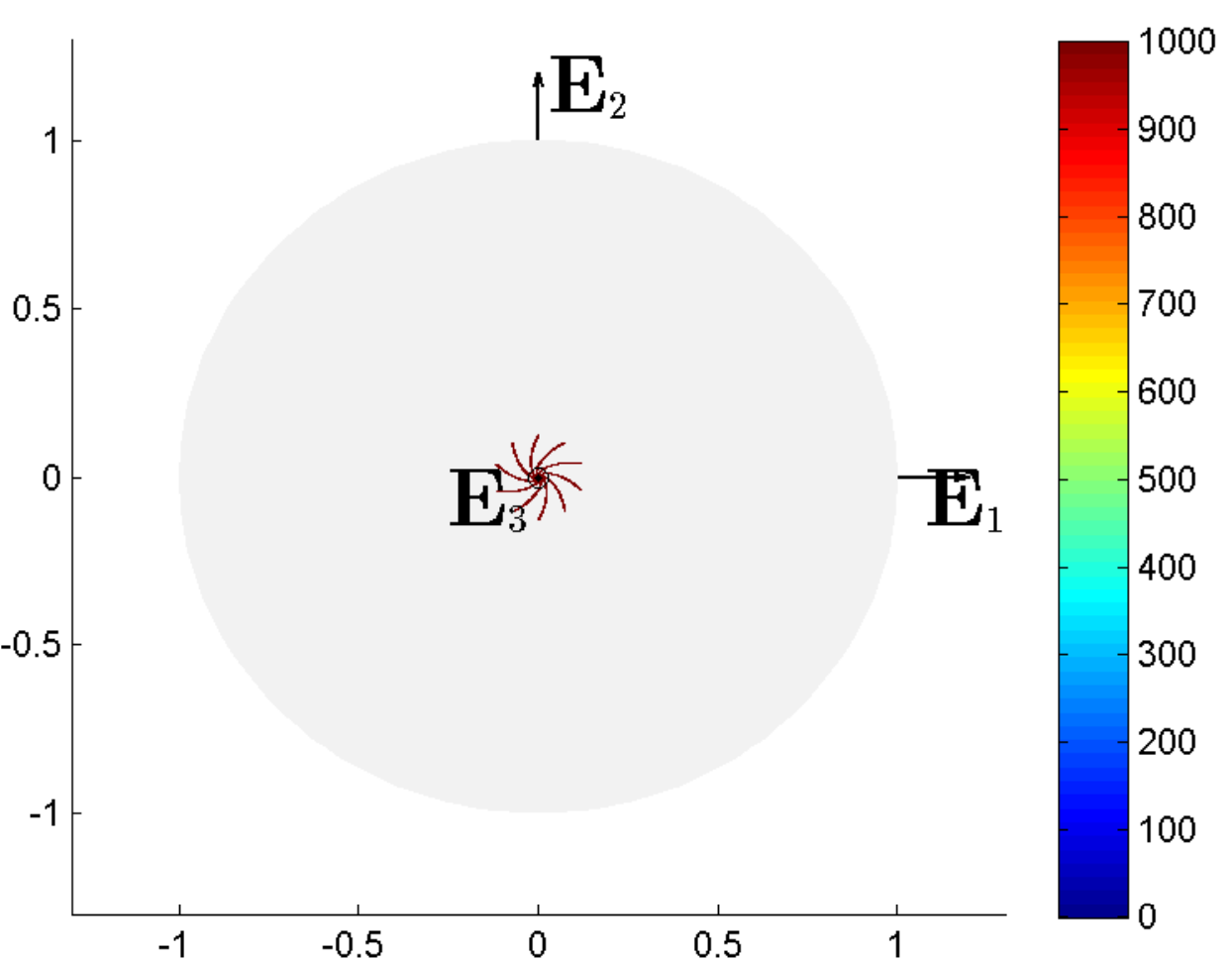}}
~
\subfloat[$t=0.0068$.\label{fig:focus_2}]{\includegraphics[width=0.49\columnwidth]{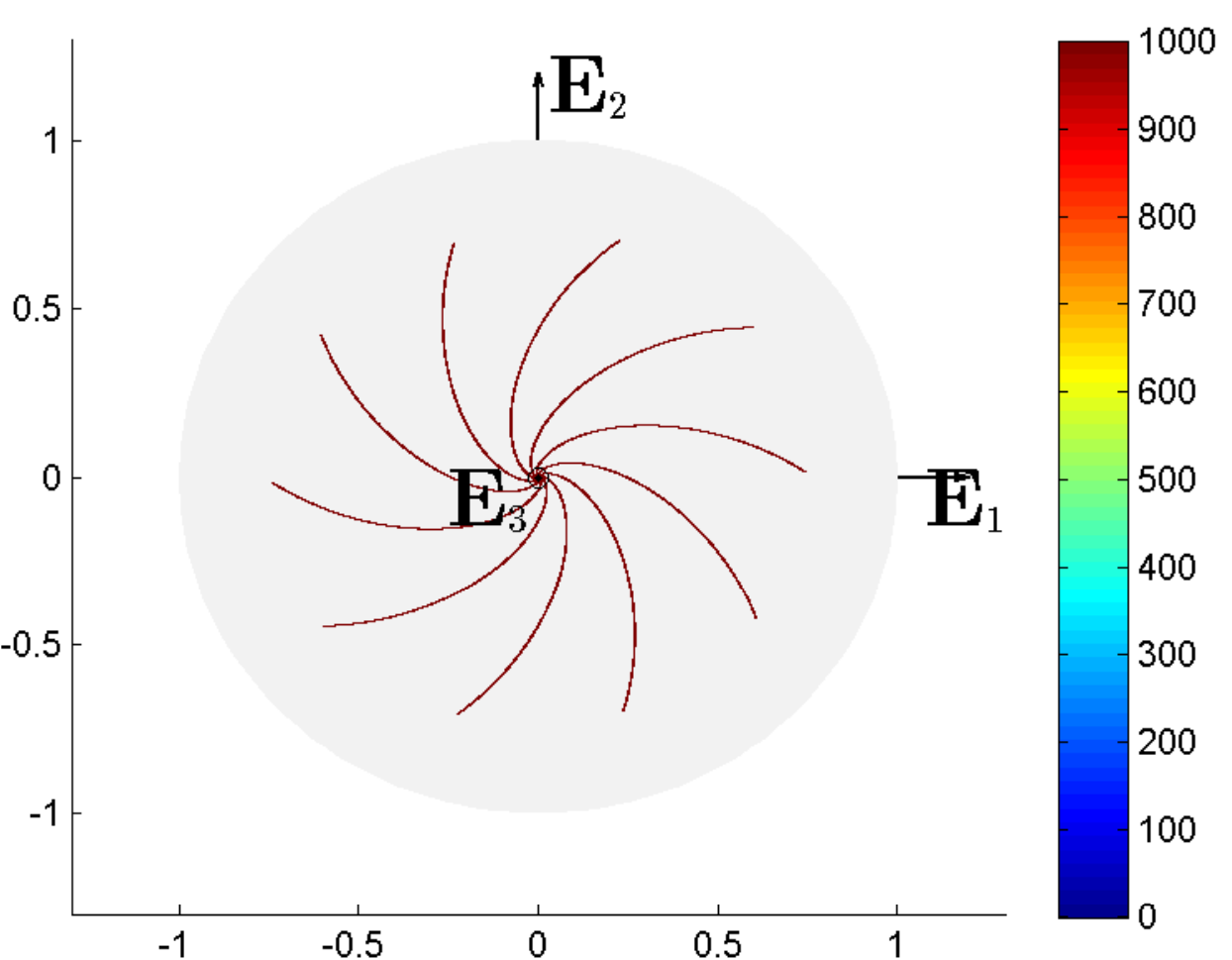}}

\subfloat[$t=0.0083$.\label{fig:focus_3}]{\includegraphics[width=0.49\columnwidth]{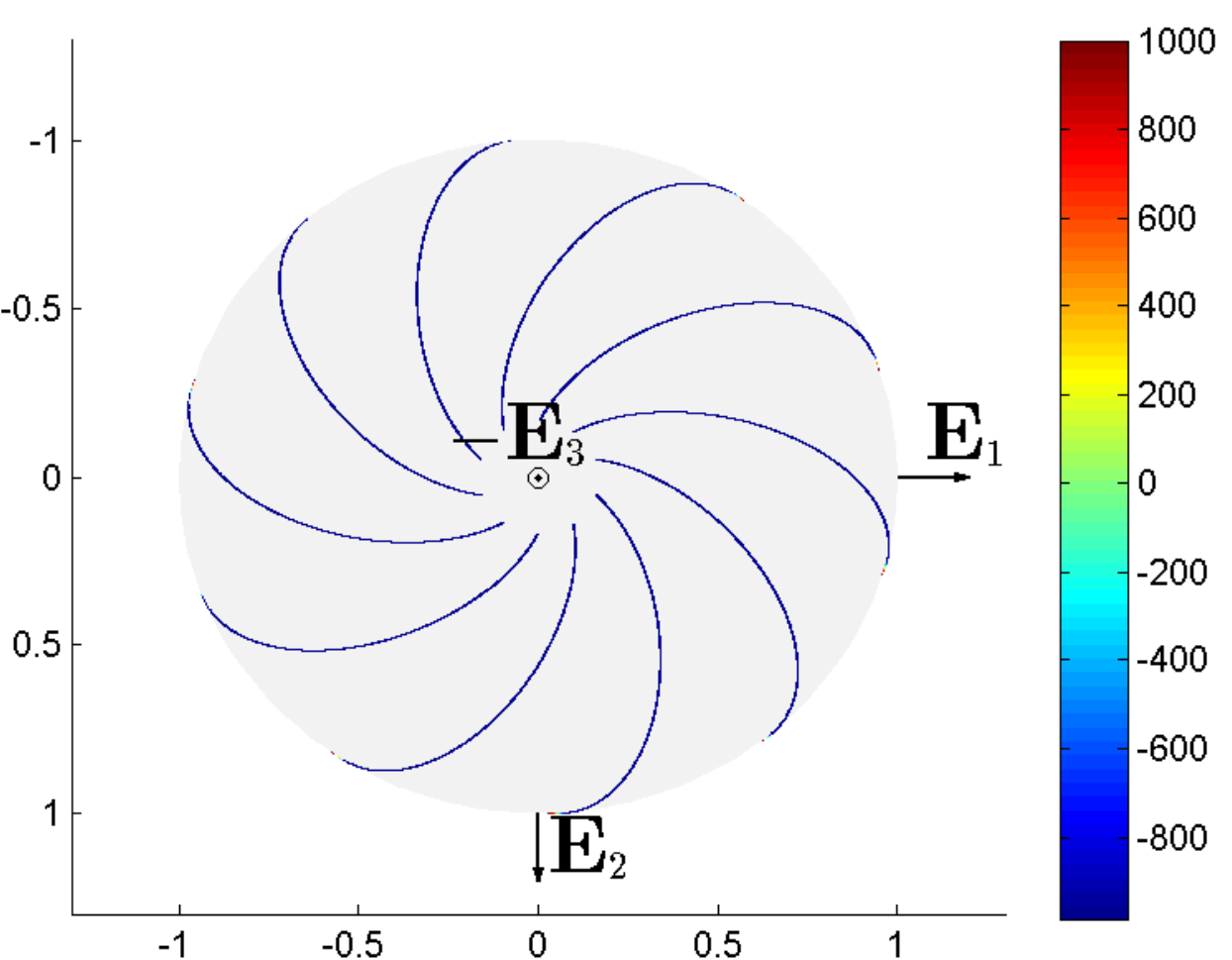}}
~
\subfloat[$t=0.0250$.\label{fig:focus_4}]{\includegraphics[width=0.49\columnwidth]{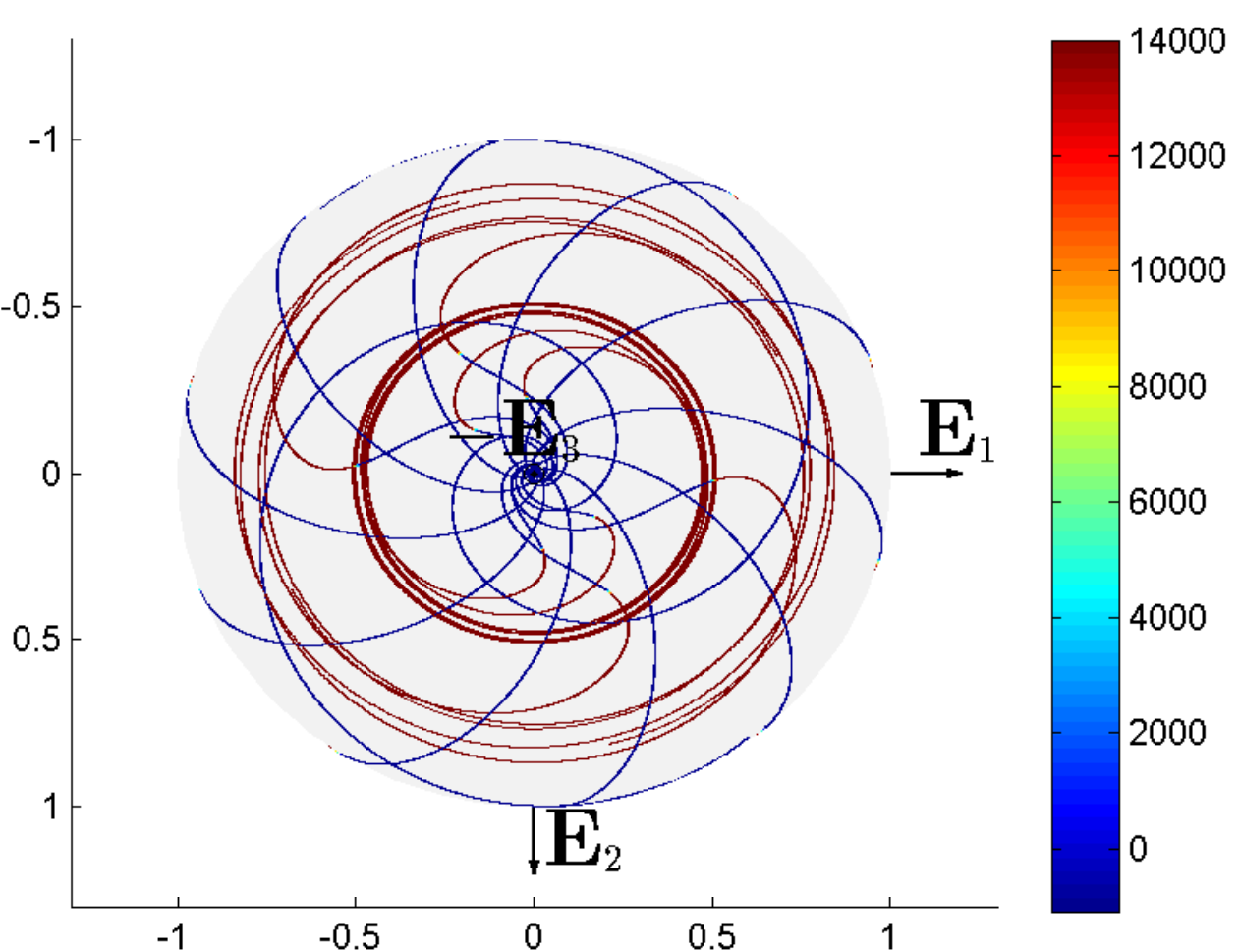}}

\subfloat[Detail near $-\mathbf{q}_{d}$\label{fig:focus_5}]{\includegraphics[width=0.49\columnwidth]{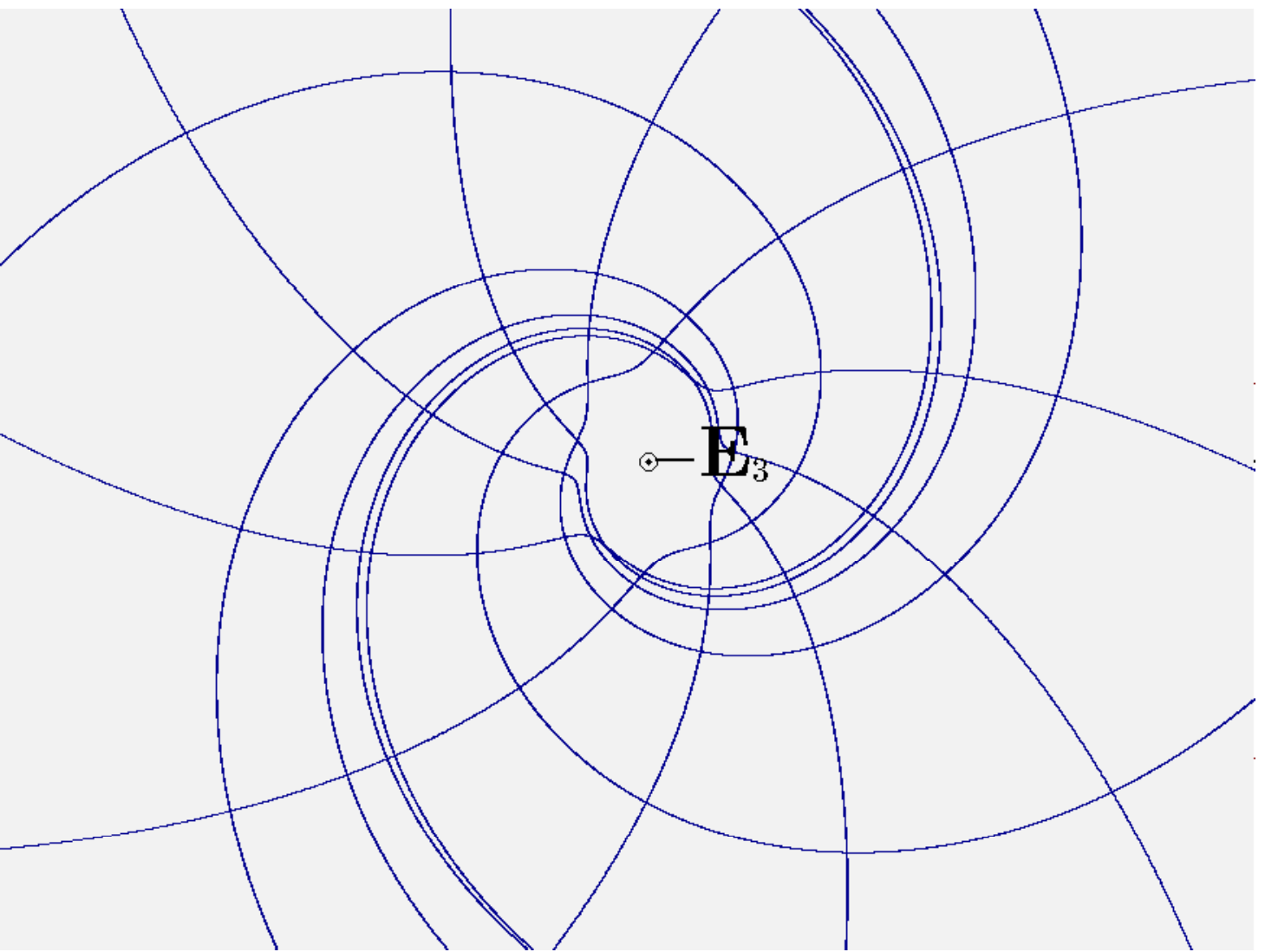}}
~
\subfloat[Distance metric, (\ref{eq:dmet})\label{fig:focus_6}]{\includegraphics[width=0.49\columnwidth]{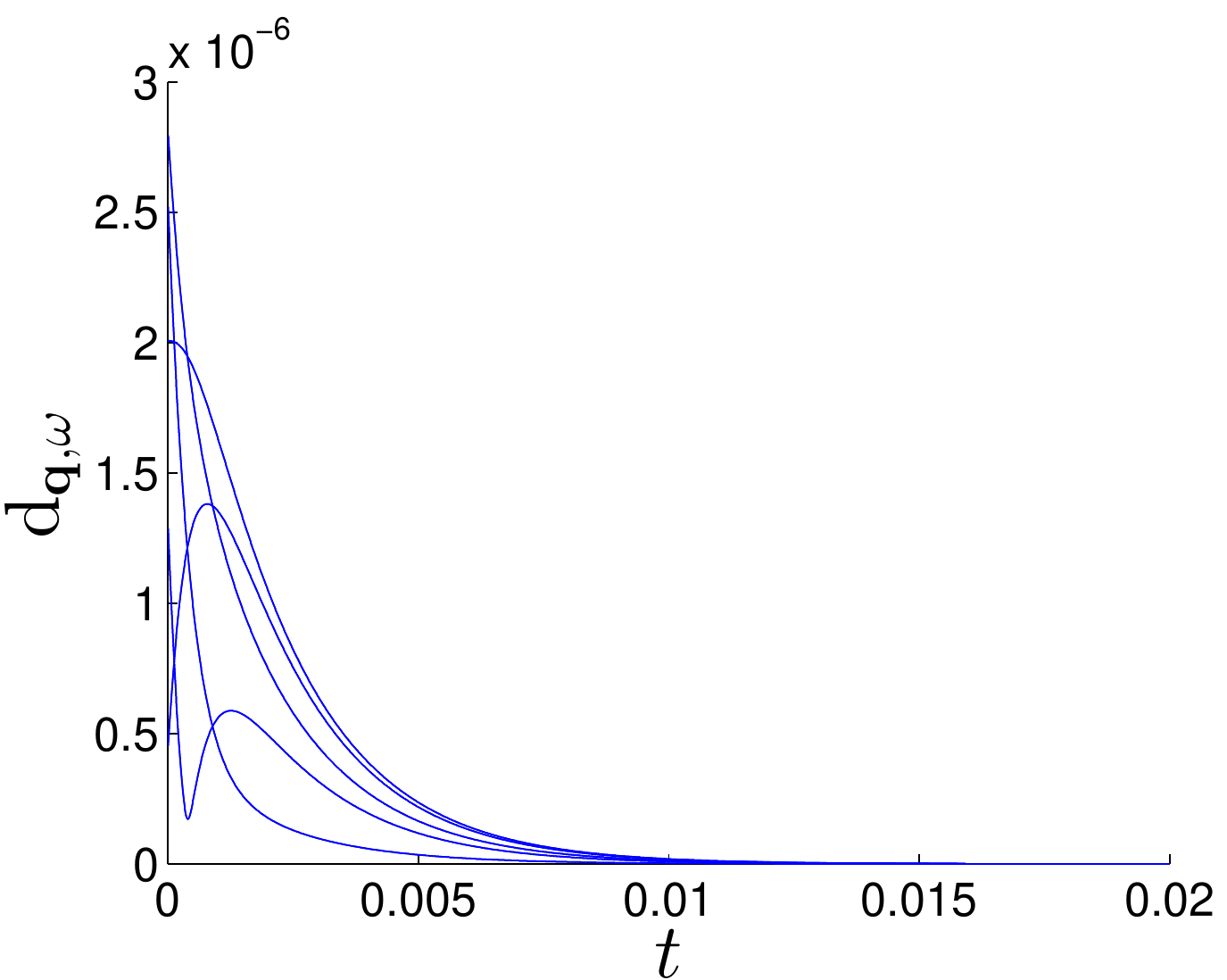}}

\caption{Backwards flow represented by $\{\boldsymbol{\varphi}^{{-}t}(\mathbf{p})\}_{t{>}0}$, $\mathbf{p}{\in}(\ref{eq:deslocs})$.
The trajectories are depicted on $\text{S}^{2}$ with ${}^{b}\boldsymbol{\omega}\cdot\mathbf{e}_{3}$ illustrated by the color of the curves according to the colorbar.
(\ref{fig:focus_1}-\ref{fig:focus_4}) Ten points from (\ref{eq:deslocs}), integrated backwards in time and shown for instances of $t$.
(\ref{fig:focus_5}) Detail of the flow near the antipodal equilibrium.
(\ref{fig:focus_6}) Forward flow represented by $\{\boldsymbol{\varphi}^{t}(\mathbf{p})\}_{t{>}0}$, $\mathbf{p}{\in}(\ref{eq:deslocs})$.
}
\label{fig:focus}
\end{figure}

On the grounds of gaining deeper understanding of the global stabilization properties of (\ref{eq:cld}), the closed-loop vector fields converging to the desired equilibrium are also visualized.
This investigation is extremely important, more important than the investigation of Section \ref{sec:fcse}, since this open dense set and the trajectories produced by the vector field dominate the dynamic behavior of the closed-loop system.
Thus an intimate knowledge of the dominant closed-loop vector field bestows to the control engineer the ability to anticipate the response of the system, a critical skill in an experimental implementation.
Because the desired equilibrium is a stable focus, we do not utilize the eigenvectors from Section \ref{sec:des}. It is sufficient to select points ''close'' enough to (\ref{eq:des}), as:
\begin{IEEEeqnarray}{L}
(\mathbf{q}_{\xi},{}^{b}\boldsymbol{\omega}_{\delta\omega})|_{(\mathbf{q}_{d},{}^{b}\boldsymbol{\omega}_{d})}{=}\Big\{(\mathbf{q},{}^{b}\boldsymbol{\omega})\in\text{S}^{2}{\times}\mathbb{R}^{3}|\IEEEnonumber\\
\mathbf{q}=\exp\left( S(\varepsilon(\cos(\vartheta)\mathbf{e}_{1}+\sin(\vartheta)\mathbf{e}_{2})\right)\mathbf{q}_{d},\IEEEnonumber\\
{}^{b}\boldsymbol{\omega}={}^{b}\boldsymbol{\omega}_{d}+\varepsilon(\cos(\vartheta)\mathbf{e}_{1}+\sin(\vartheta)\mathbf{e}_{2})+\varsigma\mathbf{e}_{3},\IEEEnonumber\\
\varepsilon=1{\cdot}10^{-6},\varsigma=1{\cdot}10^{-7},\vartheta{\in}[0,2\pi)\Big\}\label{eq:deslocs}
\end{IEEEeqnarray}
Similar to Section \ref{sec:des}, we pick ten points from (\ref{eq:deslocs}), and use the backward flow map $\boldsymbol{\varphi}^{-t}$ in evolving backwards in time the trajectories that converge to the selected points.
The generated trajectories of the flow of (\ref{eq:cld}) are shown in Fig. (\ref{fig:focus}).

Several observations regarding the trajectories of the dominant dense set are summarized next:

The trajectories are spirals (Fig. \ref{fig:focus_1}-\ref{fig:focus_3}) that as they move away from the desired equilibrium and past the antipodal equilibrium they are drawn into circular orbits that eventually wrap around $\text{S}^{2}$ as $t$ gets sufficiently large, see Fig. \ref{fig:focus_4}.
The requirement of high precision trajectory tracking demanded the use of high valued gains, see (\ref{eq:gains}).
As a result, the flow evolves extremely fast, see Fig. \ref{fig:focus_4}.
In regards to the angular velocity about the pointing direction, as $\mathbf{q}$ moves away from the desired equilibrium, ${}^{b}\boldsymbol{\omega}{\cdot}\mathbf{e}_{3}$ diverges from ${}^{b}\boldsymbol{\omega}_{d}$.
Moreover ${}^{b}\boldsymbol{\omega}{\cdot}\mathbf{e}_{3}$ changes sign as $\mathbf{q}$ passes the equator (see change in color from Fig. \ref{fig:focus_2} to Fig. \ref{fig:focus_3}), remains negative as $\mathbf{q}$ moves past the antipodal equilibrium and finally begins to increase as $\mathbf{q}$ moves in circular orbits or wraps around $\text{S}^{2}$.

The influence of the antipodal equilibrium and $\mathbf{W}^{s}(-\mathbf{q}_{d},-{}{}^{b}\boldsymbol{\omega}_{d})$ in the evolution of the solutions is apparent.
The spiral trajectories produced by the backward flow map $\boldsymbol{\varphi}^{-t}$ transform into circular orbits and take intricate paths as $t$ gets sufficiently large, see Fig. \ref{fig:focus_4}-\ref{fig:focus_5}.
More importantly, the trajectories indicate that a strong possibility exists during a step PDAV maneuver that $\mathbf{q}$ will first pass close to the antipodal equilibrium before finally converging to the desired equilibrium.
This possibility increases if $\Psi(0){>}1$ i.e., for pointing step commands of $90^{o}$ or more wrt., an equivalent axis angle rotation.
Thusly, this controller is better suited for high precision trajectory tracking, i.e., $\Psi(0){\ll}1$.

\section{CONCLUSIONS}
Using global analysis and simulation techniques, the smooth closed-loop vector fields induced by the geometric PDAV controller from \cite{VTUAV}, was visualized to gain a deeper understanding of its global stabilization properties.
A coordinate-free form of the closed-loop linearized dynamics was calculated and the local stability of each equilibrium was analyzed.
Using the solution of the linearized system, an estimate of the frequency of the complex precession/nutation oscillations that arise during PDAV trajectory tracking was obtained.
This constitutes an important tool for actuator selection and as a criterion on the feasibility of an experimental implementation.
Finally, through the use of variational integration schemes, the flow converging to the desired equilibrium and the flow ''close'' to the stable manifold of the saddle equilibrium was visualized and analyzed.
This analysis allowed the extraction of numerous observations regarding the shape of the flow, the profile of the angular velocity, the transient behavior of the solutions and finally the influence of the saddle equilibrium on the evolution of the solutions.
These observations offer intimate knowledge of the closed-loop vector field bestowing to the control engineer the ability to ''anticipate'' the response of the system, a critical skill in an experimental implementation.
Since the PDAV controller from \cite{VTUAV} can be applied to general robotic platforms, the insights and understanding gained by this analysis are applicable to a broad range of systems that utilize it and not only to a vectoring out-runner motor.






\section*{Appendix}  
Cross product map identifying the Lie algebra $\mathfrak{so}(3)$ with $\mathbb{R}^{3}$.
For $\mathbf{r}\in\mathbb{R}^3$:
\begin{IEEEeqnarray}{rCl}
\begin{array}{c}
S(\mathbf{r}){=}[{0},{-r_{3}},{r_{2}};{r_{3}},{0},{-r_{1}};{-r_{2}},{r_{1}},0]\\
{S^{-1}}(S(\mathbf{r})){=}\mathbf{r}
\end{array}\label{iso}
\end{IEEEeqnarray}
Exponential map using the Rodrigues formulation \cite{RBAC},
\begin{IEEEeqnarray}{C}
\text{exp}(\epsilon S(\boldsymbol{\xi}))=\mathbf{I}+S(\boldsymbol{\xi})\sin{\epsilon}+S(\boldsymbol{\xi})^{2}(1-\cos{\epsilon})\label{expm}
\end{IEEEeqnarray}
Time derivative of (\ref{eq:psi}) and (\ref{eq:eq}) respectively,
\begin{IEEEeqnarray}{rCl}
\dot{\Psi}&=&\mathbf{R}{}^{b}\mathbf{e}_{q}\cdot\mathbf{R}{}^{b}\mathbf{e}_{\omega}\label{eq:dpsi}\\
{}^{b}\dot{\mathbf{e}}_{q}&=&\mathbf{R}^{T}\Big(\dot{\mathbf{q}}_{d}\times\mathbf{q}+\mathbf{q}_{d}\times\dot{\mathbf{q}}\Big)-S({}^{b}\boldsymbol{\omega}){}^{b}\mathbf{e}_{q}\label{eq:deq}
\end{IEEEeqnarray}
Perturbed error vectors ${}^{b}\mathbf{e}_{q}^{\epsilon},{}^{b}\dot{\mathbf{e}}_{q}^{\epsilon},{}^{b}\mathbf{e}_{\omega}^{\epsilon}$ used during linearization,
\begin{IEEEeqnarray}{rCl}
{}^{b}\mathbf{e}_{q}^{\epsilon}&=&\mathbf{R}^{T}\Big(S(S(\mathbf{q}_{d})\mathbf{q})-S(\mathbf{q}_{d})S(\mathbf{q})\Big)\boldsymbol{\xi}\label{eq:deqe}\\
{}^{b}\dot{\mathbf{e}}_{q}^{\epsilon}&=&\mathbf{R}^{T}S(S(\dot{\mathbf{q}}_{d})\mathbf{q}+S(\mathbf{q}_{d})\dot{\mathbf{q}})\boldsymbol{\xi}
-\mathbf{R}^{T}S(\dot{\mathbf{q}}_{d})S(\mathbf{q})\boldsymbol{\xi}\IEEEnonumber\\
&&-\mathbf{R}^{T}S(\mathbf{q}_{d})S(\mathbf{R}S({}^{b}\boldsymbol{\omega})\mathbf{e}_{3})\boldsymbol{\xi}
-S(\mathbf{R}^{T}\mathbf{q}_{d})S(\mathbf{e}_{3})\delta^{b}\boldsymbol{\omega}\IEEEnonumber\\
&&+S({}^{b}\mathbf{e}_{q})\delta^{b}\boldsymbol{\omega}-S({}^{b}\boldsymbol{\omega}){}^{b}\mathbf{e}_{q}^{\epsilon}\label{eq:Ddeqe}\\
{}^{b}\mathbf{e}_{\omega}^{\epsilon}&=&\delta^{b}\boldsymbol{\omega}-\mathbf{R}^{T}S(\mathbf{R}_{d}{}^{b}\boldsymbol{\omega}_{d})\boldsymbol{\xi}\label{eq:deoe}
\end{IEEEeqnarray}

\end{document}